\documentclass{elsart}
\usepackage{amscd,amssymb,relsize,amsmath,amsfonts,latexsym}
\usepackage[arrow,matrix,graph,frame,poly,arc,tips]{xy}
\journal{Advances in Mathematics}

\newtheorem{theorem}[subsection]{Theorem}
\newtheorem{lemma}[subsection]{Lemma}
\newtheorem{proposition}[subsection]{Proposition}
\newtheorem{corollary}[subsection]{Corollary}
\newtheorem{example}[subsection]{Example}
\newtheorem{theo}{Theorem} 

\numberwithin{equation}{section}

\newcommand{\abs}[1]{\left|#1\right|}
\newcommand{\set}[1]{\left\{#1\right\}}

\renewcommand{\b}[1]{\mathbf{#1}}

\newcommand{\A}{\mathcal{A}}

\newcommand{\PP}{\mathcal{P}}
\newcommand{\RR}{\mathcal{R}}

\newcommand{\R}{\mathbb{R}}
\newcommand{\C}{\mathbb{C}}
\newcommand{\Q}{\mathbb{Q}}

\newcommand{\Z}{\mathbb{Z}}
\renewcommand{\k}{\Bbbk}

\newcommand{\CP}{{\mathbb{CP}}}

\newcommand{\sD}{{\sf D}}

\DeclareMathOperator{\rank}{rank}

\DeclareMathOperator{\Hom}{Hom}

\DeclareMathOperator{\Sym}{Sym}

\DeclareMathOperator{\Hilb}{Hilb}

\DeclareMathOperator{\nbc}{\mathbf{nbc}}
\DeclareMathOperator{\Pf}{Pf}
\DeclareMathOperator{\Pfaff}{Pfaff}
\DeclareMathOperator{\Tors}{Tors}
\DeclareMathOperator{\coker}{coker}

\DeclareMathOperator{\tc}{tc}
\DeclareMathOperator{\zcl}{zcl}
\DeclareMathOperator{\cat}{cat}
\DeclareMathOperator{\secat}{secat}
\DeclareMathOperator{\cl}{cl}
\DeclareMathOperator{\cd}{cd}

\newcommand{\tp}{{\mathsmaller{\top}}}

\begin{document}
\begin{frontmatter}

\title{Boundary manifolds of projective hypersurfaces}

\author[lsu]{Daniel C. Cohen\thanksref{nsa}}
\ead{cohen@math.lsu.edu}
\ead[url]{http://www.math.lsu.edu/\~{}cohen}
\thanks[nsa]{Partially supported 
by National Security Agency grant MDA904-00-1-0038 and 
a Faculty Research Grant from Louisiana State University}
\author[nu]{Alexander~I.~Suciu\thanksref{nsf}}
\ead{a.suciu@neu.edu}
\ead[url]{http://www.math.neu.edu/\~{}suciu}
\thanks[nsf]{Partially supported by 
National Science Foundation grant DMS-0311142 
and an RSDF grant from Northeastern University}

\address[lsu]{Department of Mathematics, Louisiana 
State University, Baton Rouge, LA 70803}
\address[nu]{Department of Mathematics, Northeastern 
University, Boston, MA 02115}


\begin{abstract}
We study the topology of the boundary manifold of a regular
neighborhood of a complex projective hypersurface.  We
show that, under certain Hodge theoretic conditions, 
the cohomology ring of the complement of the 
hypersurface functorially determines that of the boundary. 
When the hypersurface defines a hyperplane arrangement, 
the cohomology of the boundary is completely determined 
by the combinatorics of the underlying arrangement 
and the ambient dimension.  
We also study the LS category and topological complexity 
of the boundary manifold, as well as the resonance varieties 
of its cohomology ring.
\end{abstract}


\end{frontmatter}

\section{Introduction} 
\label{sec:intro}

\subsection{Boundary manifolds}
\label{intro: bdry} 

There are many ways to understand the topology of  
a homogeneous polynomial $f\colon \C^{\ell+1}\to \C$.  
The most direct approach is to study the hypersurface 
$V$ in $\CP^\ell$ defined as the zero locus of $f$. 
Another approach is to view the complement, 
$X=\CP^{\ell}\setminus V$, as the primary object 
of study.  And perhaps the most thorough is to study 
the Milnor fibration 
$f\colon \C^{\ell+1} \setminus \{f(\mathbf{x})=0\} \to \C^*$. 
Of course, the different approaches are interrelated.   
For example, if the degree of $f$ is $n$, then 
the Milnor fiber $F=f^{-1}(1)$ is a cyclic 
$n$-fold cover of $X$.  Consequently, knowledge 
of the cohomology groups of $X$ with coefficients in 
certain local systems yields the cohomology groups of $F$.  

In this paper, we take a different (yet still related) tack. 
We consider the {\em boundary manifold}, $M$, defined 
as the boundary of a closed regular neighborhood $N$ of the 
subvariety $V\subset \CP^\ell$, see Durfee~\cite{Durfee}.  
Clearly, $X\simeq  \CP^\ell \setminus N^\circ$, and 
$M$ is the boundary of  $\CP^\ell \setminus N^\circ$. 
While the complement $X$ has the homotopy type of a 
CW-complex of dimension at most $\ell$, the boundary manifold 
$M$ is a smooth, compact manifold of dimension $2\ell-1$.  

There are many questions one can ask about the topology 
of $M$, for instance, concerning its fundamental group, and 
how it relates to the fundamental group of $X$.   
In the case where $V$ is the union of an arrangement 
of lines in $\CP^2$, work in this direction was done 
by Jiang-Yau \cite{JY98}, Westlund \cite{We}, 
and Hironaka~\cite{Hi}.  Here, we resolve the 
asphericity question for the boundary manifold 
of an arbitrary hyperplane arrangement 
(see Propositions \ref{prop:nonaspherical} 
and \ref{prop:aspherical}), leaving a more detailed study of 
the fundamental group and related invariants to future work.  

For a general hypersurface $V$, our main goal in this paper 
is to compute the cohomology ring of the boundary manifold 
$M$.  We show that, under fairly mild hypotheses, the 
cohomology ring of the complement $X$ functorially 
determines the cohomology ring of $M$, and derive 
a number of consequences.  For instance, when the 
hypersurface $V=\bigcup_{H\in \A}H$ is determined 
by an arrangement of hyperplanes $\A$, these (Hodge 
theoretic) hypotheses are satisfied, and the cohomology 
of $X=X(\A)$ is thoroughly understood, thanks to 
classical results of 
Brieskorn and Orlik-Solomon.  Our results then yield 
an explicit description of the cohomology ring of the 
boundary manifold $M=M(\A)$.

\subsection{Cohomology ring of the boundary}
\label{intro: coho ring} 

Given a finite-dimensional graded algebra $A$ over 
a ring $R$, we construct a new algebra, $\sD(A)$.  
This is a particular case of a more general construction, 
the ``principle of idealization''  
due to Nagata \cite{Na}, and popularized by Reiten \cite{Re}, 
which associates to a ring $A$ and an $A$-bimodule 
$B$ the {\em trivial extension} ring $A\ltimes B:= A \oplus B$, 
with multiplication $(a,b)(a',b')=(a a',a\cdot b'+b \cdot a')$.
Applying this construction to a graded (commutative) algebra 
$A=\bigoplus_{k=0}^\ell A^k$  and the $A$-bimodule 
$B=\bar{A}=\bigoplus_{k=\ell-1}^{2\ell-1}\Hom(A^{2\ell-k-1},R)$ 
yields a graded (commutative) algebra $\sD(A)=A \ltimes \bar{A}$, 
which we refer to as the \emph{double} of $A$.

If $V \subset \CP^\ell$ is a projective hypersurface, then the 
cohomology groups of $V$ (with complex coefficients), and those 
of the complement $X=\CP^\ell \setminus V$ admit mixed Hodge 
structures.  For each $k\ge 0$, there is an increasing weight filtration 
$\{W_m\}_{m \le 2k}$ of the $k$-th cohomology group, such that 
each quotient $W_m/W_{m-1}$ has pure Hodge structure of 
weight $m$.  Our main results, proved in Section \ref{sec:ring},  
may be summarized as follows. 

\begin{theo}
\label{theo:cohomology double}
Let $V$ be a hypersurface in $\CP^\ell$, with complement $X$ and 
boundary manifold $M$.  If either $V$ is irreducible, or 
the weight filtration on the top cohomology group of $X$ 
satisfies $W_{\ell+1}(H^\ell(X;\C))=0$, then the cohomology 
ring of the boundary manifold is isomorphic to 
the double of the cohomology ring of the complement:
\begin{equation}
\label{eq:doubling}
H^*(M;\C) \cong \sD(H^*(X;\C)).
\end{equation}
\end{theo}

If $H^\ell(X;\C)$ satisfies the above weight condition and the integral 
cohomology of $X$ is torsion-free, our results can be used 
to show that the splitting \eqref{eq:doubling} 
holds over the integers, $H^*(M;\Z) \cong \sD(H^*(X;\Z))$.  
This is the case, for example, when $X$ is the complement of a 
hyperplane arrangement (see Theorem \ref{thm:coho bdry arr}).   
On the other hand, this splitting can fail with integral coefficients 
when $V$ is irreducible (see Example \ref{ex:boundary smooth}).  
With complex coefficients, the splitting \eqref{eq:doubling} can fail 
if neither of the hypotheses stated in the theorem holds (see below).  

\subsection{Arrangements and curves}
\label{intro:discuss} 
When applied to a complex hyperplane arrangement, our result 
yields an analog for the boundary manifold of a well known theorem 
of Orlik and Solomon~\cite{OS} concerning the cohomology 
ring of the complement.  Let $\A$ be a hyperplane arrangement in 
$\CP^\ell$, with complement $X(\A)$ and boundary manifold $M(\A)$.  
The integral cohomology of the complement, $H^*(X(\A);\Z)$, is torsion-free, 
and the ring structure is 
completely determined by the intersection poset $L(\A)$.  
Moreover, by work of Shapiro \cite{Sh93} and Kim \cite{Kim}, 
the complex cohomology $H^k(X(\A);\C)$ is pure of weight $2k$ for each $k$, $0\le k \le \ell$.  
It follows that the 
integral cohomology ring of the boundary manifold, $H^*(M(\A);\Z) \cong \sD(H^*(X(\A);\Z))$,  
is determined by the intersection poset $L(\A)$ and the 
ambient dimension $\ell$, see Corollary~\ref{cor:coho bdry comb}.

For an algebraic curve $V \subset \CP^2$ (in particular, 
an arrangement of lines in $\CP^2$), the associated 
boundary manifold $M$ is a Waldhausen graph manifold.  
We show in Theorem \ref{thm:bdry curve} that the ``doubling" 
formula \eqref{eq:doubling} holds for a reducible curve $V$ 
if and only if  all its components are rational curves.

Cohomology rings of graph manifolds (with $\Z_2$ coefficients)
have been the object of substantial recent study, see
Aaslepp, et.al.~\cite{ADHSZ}.  For those graph manifolds which 
arise as boundary manifolds of arrangements of rational curves in 
$\CP^2$, our methods, together with Cogolludo's computation 
of the cohomology ring of the complement of such an 
arrangement in \cite{Cog}, provide an efficient alternative.

\subsection{LS category and topological complexity}
\label{intro: ls and tc} 

Let $X^{I}$ be the space of continuous paths from the unit 
interval to $X$, and let $\pi\colon X^{I} \to X \times X$  be 
the map sending a path to its endpoints.   In \cite{Fa03}, Farber 
defines the {\em topological complexity} of $X$, denoted by $\tc(X)$, 
to be the smallest integer $k$ such that $X \times X$ can be 
covered by $k$ open sets, over each of which $\pi$ has a section.   
This numerical invariant, which depends only on the 
homotopy type of $X$, is related to the Lusternik-Schnirelmann 
category by the inequalities $\cat(X) \le \tc(X)\le 2\cat(X)-1$.  
Computing the topological complexity of $X$ is crucial to solving 
the motion planning problem for the space $X$, see \cite{Fa03}.  

The topological complexity $\tc(X)$ admits a cohomological lower 
bound in terms of the zero-divisor length of $H^*(X;\k)$, 
similar to the well known cup-length lower bound for 
$\cat(X)$.   In the case when $X=X(\A)$ is the complement 
of a hyperplane arrangement, explicit computations of $\tc(X)$ 
were carried out by Farber and Yuzvinsky \cite{FY}.  
In Section \ref{sec: top complexity}, we compute the topological 
complexity of the boundary manifold $M=M(\A)$ for various 
classes of hyperplane arrangements, using our description of the 
cohomology ring of $M$ and results from \cite{Fa03}. 
In particular, we show that the difference $\tc(M)-\cat(M)$ can 
be made arbitrarily large, see Corollary \ref{cor:tc-cat}.

\subsection{Resonance}
\label{intro: res} 

We conclude with a comparison of certain ring-theoretic invariants 
of the cohomology ring of the complement to those of the 
cohomology ring of the boundary manifold. 

Suppose $A$ is a finite-dimensional, graded, connected 
algebra over an algebraically closed field $\k$ of 
characteristic $0$.   For each $a\in A^1$, multiplication 
by $a$ defines a cochain complex $(A,a)$.  
The \emph{resonance varieties} of $A$ are the 
jumping loci for the cohomology of these complexes:
$\RR^{k}_{d}(A)=\{ a\in A^1 \mid \dim_{\k} H^k(A,a) \ge d\}$.  

In Section \ref{sec: res boundary}, we study the resonance 
varieties of the trivial extension, $\sD(A)=A\ltimes \bar{A}$. 
As an application, we obtain information about the structure 
of the resonance varieties of the boundary manifold of a 
hyperplane arrangement $\A$.  Let $A=H^*(X(\A);\k)$ be 
the Orlik-Solomon algebra.  It is well known that the components 
of the resonance varieties $\RR_d^k(X(\A))=\RR_d^k(A)$ 
are linear subspaces of $A^1 = \k^n$.  The behavior 
of the resonance varieties $\RR_d^k(M(\A))=\RR_d^k(\sD(A))$ 
is dramatically different.  Indeed, we produce examples of 
arrangements for which the resonance varieties of the 
boundary manifold contain singular, irreducible components 
of arbitrarily high degree, see Corollary \ref{cor:irr res}.

\section{The boundary manifold}
\label{sec:boundary hypersurface}

In this section, we introduce our main character, the boundary 
manifold of an (algebraic) hypersurface in complex projective space. 
We then compute its homology groups in terms of those of the 
complement to the hypersurface, and make a remark on the 
homotopy groups.

\subsection{Thickenings}
\label{subsec:thickenings}
According to C.T.C.~Wall \cite{Wa}, a thickening 
of a finite, $k$-dimen\-sional CW-complex $Y$ is 
a compact, $m$-dimensional manifold with boundary 
$W^{m}$, which is simply homotopy equivalent to $Y$.  
Such a thickening always exists, as soon as $m\ge 2k+1$:  
Embed $Y$ as a sub-polyhedron in $\R^{m}$, 
and take $W$ to be a smooth, regular neighborhood of $Y$. 

Let $M=\partial W$ be the boundary of the thickening $W$. 
In general, the homotopy type of the boundary manifold 
$M$ is not determined by the homotopy type of $Y$.  
For example, both $\CP^{2} \times D^{m-4}$ and the 
normal disk bundle of $\CP^2 \subset S^{m}$ are 
thickenings of $\CP^2$, but their boundary manifolds 
are not homotopy equivalent, see Lambrechts \cite{La00}. 
Nevertheless, if $M$ is orientable, and 
$m\ge 2(k+1)$, then the cohomology ring $H^*(M;\Z)$ is 
completely determined by $H^*(Y;\Z)$, by Poincar\'e duality 
and degree considerations. 

\subsection{Projective hypersurfaces}
\label{subsec:bdry hypersurfaces}

Let $V$ be a hypersurface in $\CP^\ell$, given as the 
zero locus of a homogeneous polynomial $f=f(\mathbf{x})$, 
where $\mathbf{x}=(x_0,\dots,x_\ell)$ are homogeneous 
coordinates on $\CP^\ell$.  A (closed) regular neighborhood, 
$N$, of $V$ in $\CP^\ell$ can be constructed either by 
triangulation, or by levels sets. In the first approach, 
triangulate $\CP^\ell$ with $V$ as a subcomplex, and 
take $N$ to be the closed star of $V$ in the second 
barycentric subdivision.  In the second, define $\phi\colon 
\CP^{\ell}\to \R$ by $\phi(\mathbf{x})=\abs{f(\mathbf{x})}^2 / 
\abs{|\mathbf{x}|}^{2d}$, where $d=\deg f$, and take 
$N=\phi^{-1}([0,\delta])$, for sufficiently small $\delta>0$. 
As shown by Durfee \cite{Durfee}, these constructions 
yield isotopic neighborhoods, independent of the choices 
made.  

Clearly, $N$ is a thickening of $V$. 
Hence, we may define {\em the boundary manifold} 
of $V$ to be
\begin{equation}
\label{eq:bdry def}
M=\partial N.
\end{equation}
This is a compact, orientable, smooth manifold of dimension $2\ell-1$. 
If $\ell=1$, then $V$ consists of, say, $n$ points on the sphere, 
and so $M$ is a disjoint union of $n$ circles.  If $\ell>1$, then 
$M$ is connected.   Here is a simple illustration.

\begin{example} 
\label{ex:boundary pencil}
{\rm 
Let $V$ be a pencil of $n+1$ hyperplanes in $\CP^{\ell}$, 
$\ell\ge 2$, defined by the polynomial $f=x_0^{n+1}-x_1^{n+1}$. 
In this case, $X$ may be realized as the complement of $n$ 
parallel hyperplanes in $\C^{\ell}$, and so it is homotopy 
equivalent to the $n$-fold wedge $\bigvee^n S^1$.  On the 
other hand, $\CP^{\ell}\setminus N=
(D^2\setminus \{\text{$n$ disks}\})\times D^{2(\ell-1)}$; 
hence $M$ is diffeomorphic to the $n$-fold connected 
sum $\#^n S^1\times S^{2(\ell-1)}$.  
}
\end{example}  

Note that the complement 
$X=\CP^\ell\setminus V$ is homotopy equivalent 
to the interior of the manifold with boundary 
$\CP^\ell \setminus N^\circ$, and that 
$M=\partial(\CP^\ell \setminus N^\circ)$.  
Also observe that, while $N$ 
is a thickening of $V$, the cohomology ring of 
$M=\partial{N}$ is 
not a priori determined by that of $V$. 

\subsection{Cohomology groups}
\label{subsec:coho hypersurfaces} 
We now analyze in detail the cohomology groups of $M$. 
We start by relating these cohomology groups  to 
those of $X$.  Throughout this section, we use integral 
coefficients, unless otherwise noted. 

\begin{proposition} 
\label{prop:coho exact seq}
Let $V$ be a hypersurface in $\CP^{\ell}$, with complement 
$X$ and boundary manifold $M$. 
For each $0\le k \le 2\ell-1$, there is an exact sequence
\begin{equation} 
\label{eq:short}
\xymatrixcolsep{16pt}
\xymatrix{0 \ar[r]& H^k(X) \ar[r]& H^k(M) \ar[r]& H^{k+1}(X,M)\ar[r]& 0}.
\end{equation}
Moreover, $H^{k+1}(X,M)\cong  H_{2\ell -k -1}(X)$, and 
the sequence splits, except possibly when $k=\ell$.
\end{proposition}

\begin{pf}
Let $i\colon M\to X$ and $j\colon V\to \CP^{\ell}$ be the inclusion 
maps.  Consider the following commuting diagram, with  
rows long exact sequences of pairs, and vertical isomorphisms 
given by the homotopy equivalence $V \hookrightarrow N$ 
and excision, respectively:
\begin{equation}
\label{eq:cd}
\xymatrixrowsep{16pt}
\xymatrixcolsep{18pt}
\xymatrix{
\ar[r] &  H^{k}(\CP^\ell,V) \ar[r] & H^k(\CP^\ell) \ar^(.55){j^*}[r] 
&  H^k(V) \ar[r] & H^{k+1}(\CP^\ell,V)  \ar[r] & 
\\
\ar[r] & H^k(\CP^\ell,N)\ar^{\cong}[d] \ar_{\cong}[u] \ar[r] 
&  H^k(\CP^\ell) \ar[d] \ar_{=}[u] \ar[r] 
& H^k(N) \ar[d] \ar_{\cong}[u] \ar[r] 
&  H^{k+1}(\CP^\ell,N)\ar^{\cong}[d] \ar_{\cong}[u] \ar[r] & 
\\
\ar[r] & H^k(X,M) \ar[r] &  H^k(X)  \ar^{i^*}[r] &  H^k(M) \ar[r] 
& H^{k+1}(X,M)  \ar[r] &\\
}
\end{equation}
By Lefschetz Duality, $H^k(\CP^\ell \setminus N^\circ,M) 
\cong H_{2\ell-k}(\CP^\ell \setminus N^\circ)$ for each $k\ge 0$.  
Since $X \simeq \CP^\ell \setminus N^\circ$, we obtain 
$H^{k}(X,M) \cong H_{2\ell-k}(X)$.  

By the Lefschetz theorem (see \cite[Ch.~5 (2.6)]{Dimca}), the map 
$j^*\colon H^k(\CP^\ell) \to H^k(V)$ is an isomorphism for 
$k\le \ell-2$ and a monomorphism for $k=\ell-1$.  
Chasing the diagram, we find that sequence 
\eqref{eq:short} is exact, for each $k\le \ell-2$. 

Now, it is well known that $X$ is a Stein space, and thus 
has the homotopy type of a CW-complex of dimension at most $\ell$.   
In particular, $H^k(X)=0$ for $k>\ell$, and $H_{\ell}(X)$ is 
finitely generated and torsion-free.  Furthermore, the boundary map 
$H^{k}(\CP^\ell,V;\Q)   \to  H^k(\CP^\ell;\Q)$ is the zero map; 
see \cite[p.~146]{Dimca}. By Lefschetz duality, 
$H^{\ell} (X,M) \cong H_{\ell}(X)$.  Hence the map   
$H^{\ell}(X,M)   \to  H^{\ell}(X)$ is the zero map. 
We conclude that sequence \eqref{eq:short} is exact 
for $k\ge \ell-1$, as well. 

For $k<\ell-1$ or $k>\ell$, one of the side terms in \eqref{eq:short} 
vanishes, so obviously the sequence splits.  For $k=\ell-1$, 
we know $H_{\ell}(X)$ is torsion-free, so \eqref{eq:short} 
splits again.  
\qed\end{pf}

\begin{corollary} 
\label{prop:cohomology groups}
The Betti numbers of the boundary manifold $M$ are given by 
$b_k(M) = b_k(X) + b_{2\ell-k-1}(X)$.  
Hence, 
the Poincar\'e polynomials of $M$ and $X$ are related by:
\begin{equation}
\label{eq:poincare polynomial}
P(M,t) = P(X,t) + t^{2\ell-1}\cdot P(X,t^{-1}).
\end{equation}
\end{corollary}

Proposition \eqref{prop:coho exact seq} determines the 
cohomology groups of $M$ in terms of the (co)ho\-mo\-logy 
groups of $X$, except possibly the torsion in $H^\ell(M)$.  
By the Universal Coefficient Theorem, this torsion 
fits into the short exact sequence
\begin{equation}
\label{eq:tors}
\xymatrixcolsep{16pt}
\xymatrix{0  \ar[r] &\Tors(H_{\ell-1}(X)) \ar[r] & \Tors(H^{\ell}(M)) \ar[r] 
& \Tors(H_{\ell-1}(X))  \ar[r] & 0}.
\end{equation}
This sequence may or may not split.  As we shall see in 
examples below, both possibilities can occur.

\begin{example}
\label{ex:boundary smooth}
{\rm
Let $V$ be a smooth algebraic hypersurface in $\CP^{\ell}$ 
of degree $d$.  In this case, $N$ can be taken to be 
a tubular neighborhood of $V$, diffeomorphic to the 
unit normal disk bundle $\nu$.  Hence $M$ is the total 
space of the $S^1$-bundle over $V$ with Euler number
$e=c_1(\nu)[V]$.  

In particular, if $\ell=2$, then $V$ is a curve of genus 
$g=\binom{d-1}{2}$, with $e=d^2$.  Hence, by the Gysin 
sequence, $H^2(M)=\Z_{d^2}$.  On the other hand, 
$H_1(X)=\Z_d$. Thus, in this instance, \eqref{eq:tors} 
is a non-split exact sequence, of the form 
$0  \to \Z_d\to \Z_{d^2} \to\Z_d  \to 0$.
}
\end{example}

\subsection{Affine hypersurfaces and Milnor fibrations}
\label{subsec: affine and milnor}

Let $V_0 \subset \C^\ell$ be an affine hypersurface, 
defined by the vanishing of a polynomial 
$f_0=f_0(x_1,\dots,x_\ell)$ of degree $n$. 
Let $V$ be the projective closure of $V_0$, 
defined by the vanishing of the homogeneous polynomial 
$f(x_0,x_1,\dots,x_\ell)= x_0^{n+1}\cdot f_0(x_1/x_0,\dots,x_\ell/x_0)$.  
Clearly, $\CP^\ell \setminus V=\C^\ell\setminus V_0$.  

If $f_0$ itself is homogeneous,   
then $f(x_0,x_1,\dots,x_\ell)= x_0\cdot f_0(x_1,\dots,x_\ell)$.
Moreover, we can take the regular 
neighborhood $N$ of $V$ to be the union of a regular 
neighborhood of $V_0$, say $N_0$, with a tubular neighborhood 
of the hyperplane at infinity (after rounding corners).  
Thus, $\CP^\ell \setminus N^\circ$ is diffeomorphic to 
$D^{2\ell} \setminus (D^{2\ell}\cap N_0^\circ)$, and so
\begin{equation*}
\label{eq:M decomposition}
M= (S^{2\ell -1} \setminus (S^{2\ell -1} \cap N_0))\bigcup 
D^{2\ell} \cap \partial N_0 \, .
\end{equation*}

As shown in \cite{Mil}, each of the two sides of the above 
decomposition is diffeomorphic to the total space of the 
Milnor fibration, $F\to Y\to S^1$, determined by the 
homogeneous polynomial $f_0$. Thus, $M$ is the 
double of the manifold with boundary $Y$:
\begin{equation}
\label{eq:boundary as double}
M=\partial( Y \times I) = Y \cup_{\partial Y} Y.
\end{equation}
Furthermore, $M$ fibers over the circle, with fiber the 
double of $F$. 

Notice that, in this situation, the exact sequence \eqref{eq:short} 
always splits.  Indeed, the inclusion $Y\to X$ is a homotopy 
equivalence, which factors through the inclusions $Y\to M$ 
and $i\colon M\to X$.  Thus, $i^*\colon H^*(X)\to H^*(M)$ is a 
split injection. 

\begin{example} 
\label{ex:boolean}
{\rm
Let $f=x_0x_1\cdots x_{\ell}$ be the polynomial defining 
the Boolean arrangement in $\CP^{\ell}$. 
Then $M=S^{\ell-1}\times T^{\ell}$, where $T^\ell$ 
is the $\ell$-torus; see \cite[Example~2.29]{Dimca}.
}
\end{example}  

\begin{example} 
\label{ex:boundary near pencil}
{\rm
Let $f=x_0(x_1^n-x_2^n)$ be the polynomial defining 
a near pencil of $n+1$ lines in $\CP^2$. In this case, 
$Y$ admits a fibration over the circle (different from the 
Milnor fibration!), with fiber $D^2 \setminus \{\text{$n-1$ disks}\}$, 
and monodromy a Dehn twist about the boundary $D^2$. 
It follows that $M=\Sigma_{n-1} \times S^1$, where 
$\Sigma_g$ denotes a surface of genus $g$.  
}
\end{example}  

\begin{example}
\label{ex:boundary product pencils}
{\rm
More generally, let 
$f=x_0(x_1^{n_1} - y_1^{n_1})\cdots 
(x_k^{n_k} - y_k^{n_k})$, with  $n_i\ge 2$.  Then  
$M=T^k \times (\#^{m} T^k \times S^{2k-1})$, where 
$m=\prod_{i=1}^{k} (n_i-1)$.  
}
\end{example}  

\begin{example} 
\label{ex:sum of squares}
{\rm
Let $f=x_0 (x_1^2+\cdots +x_{\ell}^2)$.  In this case, 
the Milnor fiber $F$ of $f_0=x_1^2+\cdots+x_{\ell}^2$ 
is diffeomorphic to the unit disk bundle of $S^{\ell-1}$.  
Thus, $M$ fibers over $S^1$ with fiber $E$, where 
$S^{\ell-1}\to E \to S^{\ell-1}$ is the bundle with Euler 
number $1-(-1)^{\ell}$. 

Now assume $\ell$ is odd and $\ell>1$. A computation 
with the Wang sequence for the bundle $F\to Y\to S^1$ 
shows that $H_{\ell-1}(X)=\Z_2$; see \cite[Example~3.2]{DN}. 
Hence,  \eqref{eq:tors} is a split exact sequence, of the form 
$0  \to \Z_2\to \Z_{2}\oplus\Z_2 \to\Z_2  \to 0$.
}
\end{example}  

\subsection{On asphericity of the boundary} 
\label{subsec:kpi1} 
If $V$ is a hypersurface in $\CP^\ell$, the inclusion map 
$M\to X$ is an $(\ell-1)$-equivalence, see for 
instance \cite[Proposition 2.31]{Dimca}; in particular, 
$\pi_i(M)\cong\pi_i(X)$, for $i<\ell-1$.  A natural question 
arises:  Is $M$ aspherical?  In other words, do all the higher homotopy 
groups of $M$ vanish?

If $\ell=2$, the manifold $M^3$ is a graph manifold in the sense 
of Waldhausen.  With a few exceptions (such as lens spaces), 
manifolds of this type are aspherical.  In higher dimensions, though, 
this never happens.

\begin{proposition}
\label{prop:nonaspherical}
Let $M$ be the boundary manifold of a hypersurface in $\CP^{\ell}$.  
If $\ell\ge 3$, then $M$ is not aspherical.
\end{proposition}
\begin{pf}
Let $\pi=\pi_1(M)$ be the fundamental group of $M$.  
Since the inclusion $i\colon M\to X$ is an 
$(\ell-1)$-equivalence, and since $\ell\ge 3$, the induced map 
$i_*\colon \pi_1(M)\to \pi_1(X)$ is an isomorphism.  
Let $g\colon X \to K(\pi,1)$ be a classifying map for the universal 
cover $\widetilde{X}\to X$.  
By definition, $g_*\colon \pi_1(X)\to \pi$ is an isomorphism.  
Hence, the composite map $g\circ i\colon M \to K(\pi,1)$ is a 
classifying map for $\widetilde{M}\to M$.

Now suppose $M$ is aspherical.  Then the map $g\circ i\colon M \to K(\pi,1)$ 
must be a homotopy equivalence, since it induces an isomorphism 
on fundamental groups.  Consequently, 
$(g\circ i)^*\colon H^{2\ell-1}(\pi) \to H^{2\ell-1}(M)=\Z$ is an isomorphism.  
On the other hand, $i^*\colon H^{2\ell-1}(X) \to H^{2\ell-1}(M)$ is the 
zero map, since the CW-complex $X$ has dimension at most $\ell$.  
This contradiction finishes the proof.
\qed\end{pf}

\section{The cohomology ring of the boundary manifold}
\label{sec:ring}

Let $V \subset \CP^\ell$ be a projective hypersurface, 
with complement $X=\CP^\ell\setminus V$, and associated 
boundary manifold $M$.  In this section, we determine the 
structure of the cohomology ring $H^*(M;\C)$ under certain 
conditions.  These conditions are given below in terms of the mixed 
Hodge structure on $H^*(X;\C)$, respectively $H^*(V;\C)$.  
First, we discuss the relevant algebraic structure. 

\subsection{The double of a graded ring}
\label{subsec:double}
If $A$ is a ring and $B$ is an $A$-bimodule, the 
\emph{trivial extension} of $A$ by $B$, written $A\ltimes B$, 
is the additive group $A\oplus B$, with multiplication 
given by 
$(a,b)(a',b')=(aa',a\cdot b'+b\cdot a')$, see \cite{Na,Re}.  
Note that $A \cong\{(a,0)\}$ is a subring of the trivial 
extension, and that $B\cong\{(0,b)\}$ is a square-zero ideal.

Now let $A=\bigoplus_{k=0}^\ell A^k$ be a finite-dimensional 
graded ring over a base ring $R$.  We will assume $R$ is a 
commutative ring with $1$, and all the graded pieces $A^k$ 
are finitely generated free $R$-modules.  Define the 
\emph{double} $\sD(A)$ of $A$ to be the trivial extension 
of $A$ by the graded $A$-bimodule 
$\bar{A}=\bigoplus_{k=\ell-1}^{2\ell-1} \bar{A}^k$, where 
$\bar{A}^k = \Hom(A^{2\ell-k-1},R)$, and 
the $A$-bimodule structure is given by $a\cdot b(x)=b(xa)$ 
and $b \cdot a(x)=b(ax)$ for $a,x\in A$ and $b\in \bar{A}$.  
If $A$ is a graded commutative ring, it is readily checked 
that $\sD(A)=A\ltimes \bar{A}$ is a graded commutative 
ring as well.

Let $\mu\colon A \otimes A \to A$, $\mu(a,a')=aa'$, denote the 
multiplication map of the ring $A$. Then the multiplication map 
$\sD(\mu)\colon \sD(A) \otimes \sD(A) \to \sD(A)$ of the double 
restricts to $\mu$ on $A \otimes A$ and vanishes on 
$\bar{A}\otimes \bar{A}$, while on $A\otimes \bar{A}$ 
it vanishes, except for 
\begin{equation}
\label{eq:multiplication}
\sD(\mu)(a^k_j ,\bar a^r_p) =
\sum_i \mu_{i,j,p}\, \bar{a}^{r-k}_i, \quad \text{if \ 
$\mu(a^{r-k}_i, a^k_j) =
\sum_p \mu_{i,j,p}\, a^{r}_p$},
\end{equation}
where $\{a^k_j\}$ is a (fixed) homogeneous basis for $A^k$  
and $\{\bar{a}^k_j\}$ is the dual basis for 
$\bar{A}^{2\ell-k-1}=\Hom(A^{k},\C)$. 
The proof of the next result is straightforward.

\begin{proposition} 
\label{prop:functorial}  
The doubling construction is functorial. 
In particular, if $A_1$ and $A_2$ are isomorphic as graded rings, then 
$\sD(A_1)$ and $\sD(A_2)$ are isomorphic as graded rings.
\end{proposition}

Denote the Betti numbers of $A$ by $b_k(A)=\rank A^{k}$, and let
$\Hilb (A,t)=\sum_{k=0}^{\ell} b_k(A)\cdot t^k$ 
be the Hilbert series of $A$.  Then:
\begin{equation}
\label{eq:hilb}
\Hilb(\sD(A),t) = \Hilb(A,t) + t^{2\ell-1}\cdot \Hilb(A,t^{-1}). 
\end{equation}
In particular, if $A$ is connected (i.e., $b_0(A)=1$), then $\sD(A)$ 
is an Artin-Gorenstein ring. 

Recall from Proposition \ref{prop:cohomology groups} that the (integral) 
cohomology of the boundary manifold $M$ is additively given by 
$H^q(M) \cong H^q(X) \oplus H^{q+1}(X,M)$.  Let $R$ be a 
coefficient ring. 

\begin{theorem} 
\label{thm:square-zero} 
Assume that $H^*(X;R)$ is a free $R$-module.  
If $H^*(X,M;R)$ is a square-zero subring of $H^*(M;R)$, then 
$H^*(M;R) \cong \sD(H^*(X;R))$ as graded rings.
\end{theorem}

\begin{pf}
Recall that the inclusion $i\colon M \to X$ induces 
a monomorphism 
$i^*\colon H^*(X)\to H^*(M)$ in cohomology.  Let $A=i^*(H^*(X;R))$, 
and note that $A$ is a subring of $H^*(M;R)$.  Comparing formulas 
\eqref{eq:poincare polynomial} and \eqref{eq:hilb}, and using the 
$R$-freeness assumption for $H^*(X;R)$, 
we see that $H^*(M;R)$ and $\sD(A)=A\ltimes \bar{A}$ 
are additively isomorphic.  So it suffices 
to show that the cup-product structure in $H^*(M;R)$ 
coincides with the multiplicative structure in $\sD(A)$.  
This is clearly the case for the restriction to the common 
subring $A$.

For simplicity, let us suppress the coefficient ring $R$ from the 
notation.  Fix a generator $\omega \in H^{2\ell-1}(M)$, and note 
that $\omega \notin A$.  For each $q$, $0 \le q \le \ell$, let 
$\{a^q_1,\dots,a^q_{b_q}\}$ be a basis for $A_q \cong H^q(X)$, 
where $b_q=b_q(A)$.  By Poincar\'e 
duality, there are linearly independent elements 
$\bar{a}^{q}_1,\dots,\bar{a}^{q}_{b_q}$ in $H^{\bar{q}}(M)$ 
so that $a^q_i \cup \bar{a}^{q}_j = \delta_{i,j} \omega$, where 
$\bar{q}=2\ell-q-1$ and $\delta_{i,j}$ is the Kronecker index.  
Since $A$ is a subring of $H^*(M)$ and $\omega \notin A$, 
the dual classes $\bar{a}^{q}_i$ are also not in $A$.  
Identifying $H^q(M)=H^q(X)\oplus H^{q+1}(X,M)$, it 
follows that $\{\bar{a}^q_1,\dots,\bar{a}^q_{b_q}\}$ 
forms a basis for $H^{\bar{q}+1}(X,M) \subset H^{\bar{q}}(M)$.  
Consequently, $H^q(M)$ has basis $\{a^q_1,\dots,a^q_{b_q}, 
\bar{a}^{\bar{q}}_1,\dots,\bar{a}^{\bar{q}}_{b_{\bar{q}}}\}$.

By hypothesis, we have $\bar{a}_i^p \cup \bar{a}^q_j=0$ for all 
$p,q$ and $i,j$.  It remains to consider the cup-product 
$a^p_j \cup\bar{a}^q_k \in H^{p+\bar{q}}(M)$.  If $p=0$, then 
$a^p_j \cup\bar{a}^q_k=1\cup\bar{a}^q_k=\bar{a}^q_k$.  If $p>q$, 
then $a^p_j \cup\bar{a}^q_k=0$.  So assume that $0<p \le q \le \ell$, 
which implies that $p+\bar{q} \ge \ell$.

If $p+\bar{q} > \ell$, then $a^p_j
\cup \bar{a}^q_k = \sum_{i = 1}^{b_{q-p}} c_{i,j,k}\bar{a}^{q-p}_i$
for some constants $c_{i,j,k}$.  Write the multiplication in
$A\cong H^*(X)$ as $a^r_i \cdot a^p_j = \sum_{l=1}^{b_{r+p}}
\mu_{i,j,l} a^{r+p}_l$, and note that
$\mu_{j,i,l}=(-1)^{rp}\mu_{i,j,l}$ in this instance. For a fixed
$i$, cupping with $a^{q-p}_i$ yields $a^{q-p}_i \cup
a^p_j \cup \bar{a}^q_k = c_{i,j,k}\omega$. Since
\[
a^{q-p}_i \cup a^p_j \cup \bar{a}^q_k = \Bigl(\sum_{l=1}^{b_{q-p}}
\mu_{i,j,l} a^q_l\Bigr) \cup \bar{a}^q_k = \mu_{i,j,k}\omega,
\]
we must have $c_{i,j,k}=\mu_{i,j,k}$, and so $a^p_j \cup \bar{a}^q_k =
\sum_{i=1}^{b_{q-p}} \mu_{i,j,k}\bar{a}^{q-p}_i$. 

We are left with the case $p+\bar{q} = \ell$, that is, 
$p=1$ and $q=\ell$.  We then have 
$a^1_j \cup \bar{a}^\ell_k = \sum_{i=1}^{
b_{\ell-1}} c_{i,j,k} \bar{a}^{\ell-1}_i+ \sum_{i=1}^{b_\ell}
d_{i,j,k} a^\ell_i$ for some constants $c_{i,j,k}$ and
$d_{i,j,k}$.  Since $0=\bar{a}^\ell_i \cup a^1_j \cup
\bar{a}^\ell_k = \pm d_{i,j,k} \omega$, we have $d_{i,j,k}=0$.
Then, a calculation as above yields $c_{i,j,k}=\mu_{i,j,k}$, where
$a^{\ell-1}_i \cdot a^1_j = \sum_{k=1}^{b_\ell} \mu_{i,j,k}
a^\ell_k$. Thus, $a^1_j \cup \bar{a}^\ell_k = 
\sum_{i=1}^{b_{\ell-1}} \mu_{i,j,k} \bar{a}^{\ell-1}_i$.

Notice that these calculations show that the square-zero subring 
$H^*(X,M)$ is, in fact, an ideal in $H^*(M)$.  
Using these calculations, and formula \eqref{eq:multiplication}, 
it is readily checked that the cup-product structure in $H^*(M)$ 
coincides with the multiplicative structure in $\sD(H^*(X))$.
\qed\end{pf}

The freeness assumption from Theorem~\ref{thm:square-zero} 
holds, for example, when $R=\Z$ and $H^*(X)$ is torsion-free, 
or when $R=\k$ is a field.  This assumption is necessary, as 
illustrated by the smooth plane curve of degree $d>1$ from 
Example \ref{ex:boundary smooth}.  Indeed, for such a curve, 
$H^2(M;\Z)=\Z_{d^2}$ does not split as a direct sum, and so 
$H^*(M;\Z)\not\cong \sD(H^*(X;\Z))$, even though $H^*(X,M;\Z)$ 
is a square-zero subring of $H^*(M;\Z)$, by degree considerations. 

\subsection{Hodge structures} 
\label{subsec:hodge}
Now we pursue conditions which insure that the hypotheses of 
Theorem \ref{thm:square-zero} hold.  These conditions will be 
given in terms of mixed Hodge structures.  For the rest of this section, 
we shall take coefficients in the ring $R=\C$. 

If $V$ is a smooth projective variety, then, by a classical theorem 
of Hodge, each cohomology group $H^m(V)$ admits a pure Hodge 
structure of weight $m$.  That is, for $H=H^m(V)$, there is a direct 
sum decomposition 
\begin{equation} 
\label{eq:pure weight m}
H = \bigoplus_{p+q=m} H^{p,q},
\end{equation}
where $\overline{H^{p,q}}=H^{q,p}$ (complex conjugation). 

If $X$ is a quasi-projective variety, then, by a well known 
theorem of Deligne~\cite{Del}, each cohomology group of $X$ 
admits a mixed Hodge structure.  That is, 
for each $k$, there is an increasing weight filtration
\begin{equation} \label{eq:weight filtration}
0 = W_{-1} \subset W_0 \subset \dots \subset W_{2k} = H^k(X), 
\end{equation}
such that each quotient $W_m/W_{m-1}$ of the subspaces 
$W_m=W_m(H^k(X))$ of $H^k(X)$ admits a pure Hodge 
structure of weight $m$ as in \eqref{eq:pure weight m}.

The following properties of the weight filtration will be of use.  See 
\cite{Dimca,Dur83,PeSt} for further details.
\begin{enumerate}
\item \label{item:proj} If $X$ is projective, then $W_k=H^k(X)$ 
for each $k$.
\item \label{item:smooth}
If $X$ is smooth, then $0=W_{k-1} \subset H^k(X)$ for each $k$.  
\item \label{item:compactification} For any smooth compactification 
$\iota\colon X \to \bar{X}$ of $X$, 
$W_k=\iota^*(H^k(\bar{X}))$ for each~$k$.
\item \label{item:functorial} The weight filtration is functorial.  
For an algebraic map $f\colon X\to Y$, 
the induced homomorphism $f^*$ strictly preserves the filtration:  
If $x\in W_m(H^k(X))$ is in the image of $f^*$, there is an element 
$y \in W_m(H^k(Y))$ with $f^*(y)=x$.
\end{enumerate}

It follows from work of Durfee and Hain \cite{DH88} that the 
cohomology of the boundary manifold $M$ of a projective 
hypersurface $V$ admits a mixed Hodge structure.  
Furthermore, the cup-product of $H^*(M)$ is a morphism of mixed Hodge 
structures, and the top cohomology $H^{2\ell-1}(M)$ is of weight $2\ell$ 
(and type $(\ell,\ell)$).

\begin{theorem} 
\label{thm:irr double}
Let $V$ be a hypersurface in $\CP^\ell$ with complement $X$  
and boundary manifold $M$.  If $V$ is irreducible, then 
$H^*(M;\C) \cong \sD(H^*(X);\C)$ as graded algebras.
\end{theorem}
\begin{pf}
If $\ell=1$, then $V$ is a point in $\CP^1$.  In this 
instance, $X$ is contractible, $M$ is a circle, 
and it is readily checked that $H^*(M) \cong \sD(H^*(X))$.  

So we may assume that $\ell \ge 2$.  
By Theorem \ref{thm:square-zero}, it suffices to show that $H^*(X,M)$ is a 
square-zero subalgebra of $H^*(M)$.  For this, it is enough to show that 
$u \cup v=0$ for $u\in H^{r+1}(X,M) \subset H^r(M)$ and 
$v\in H^{s+1}(X,M) \subset H^s(M)$, where $\ell-1 \le r,s \le \ell$.

Recall that, for $k\le 2\ell-2$, the inclusion $j\colon V\to \CP^\ell$ 
induces a monomorphism in $k$-th cohomology.  From diagram \eqref{eq:cd}, 
we see that $H^{k+1}(X,M)$ is isomorphic to $H^k_0(V)$, the primitive 
cohomology of $V$, given by $H^k_0(V) = \coker [j^*\colon H^k(\CP^\ell) \to 
H^k(V)]$.

It is known that the connecting homomorphism in the long 
exact sequence of the pair is weight-preserving, see \cite{Dimca,Dur83,PeSt}.  
This fact, and the properties recorded above, imply that all 
cohomology classes in $H^{k+1}(X,M)\cong 
H^k_0(V)$ (for $k\le 2\ell-2)$ are of weight at most $k$.  

Now take $u\in H^{r+1}(X,M) \subset H^r(M)$ and $v\in H^{s+1}(X,M) \subset 
H^s(M)$ as above.  If $r=s=\ell$, then clearly $u \cup v=0$.  If, say, 
$r=\ell-1$ and $s=\ell$, then $u$ is of weight at most $\ell-1$ 
and $v$ is of weight at most $\ell$.  Hence, $u\cup v$ is of 
weight at most $2\ell-1$ in $H^{2\ell-1}(M)$.  But 
$W_{2\ell-1}(H^{2\ell-1}(M))=0$ by the results of Durfee 
and Hain noted above.  So we must have $u\cup v=0$.

Finally, if $r=s=\ell-1$, then $u\cup v$ is of weight at most $2\ell-2$ in 
$H^{2\ell-2}(M)$.  Since $V$ is irreducible, 
$H^1(X)=0$, the map $j^*\colon H^{2\ell-2}(\CP^\ell)\to H^{2\ell-2}(V)$ is an 
isomorphism, and $H^{2\ell-1}(X,M)\cong H^{2\ell-2}_0(V)=0$.  If $\ell=2$, 
then all non-trivial classes in $H^2(M)=H^2(X)$ are of weight at least $3$ by 
Poincar\'e duality, since all classes in $H^1(M)\cong H^1_0(V)=H^1(V)$ are of 
weight at most $1$.  If $\ell\ge 3$, then $H^{2\ell-2}(M)=H^{2\ell-2}(X)=0$ 
since $X$ has the homotopy type of an $\ell$-dimensional complex.  
It follows that $u\cup v=0$ in either case.
\qed\end{pf}

\begin{theorem} 
\label{thm:weight double}
Let $V$ be a hypersurface in $\CP^\ell$ with complement $X$ and 
boundary manifold $M$.  
If $W_{\ell+1}(H^\ell(X;\C))=0$, then 
$H^*(M;\C) \cong \sD(H^*(X;\C))$ as graded algebras.
\end{theorem}
\begin{pf}
If $\ell=1$, then $V$ is a union of, say, $n+1$ points in $\CP^1$.  In this 
instance, $X$ is homotopic to a bouquet of $n$ circles, $M$ is a disjoint 
union of $n+1$ circles, and it is readily checked that $H^*(M) \cong 
\sD(H^*(X))$.  

If $\ell\ge 2$, by Theorem \ref{thm:square-zero}, it suffices to 
show that $H^*(X,M)$ is a square-zero subalgebra of $H^*(M)$.  
For this, as above, it is enough to show that 
$u \cup v=0$ for $u\in H^{r+1}(X,M) \subset H^r(M)$ and 
$v\in H^{s+1}(X,M) \subset H^s(M)$, where $(r,s)=(\ell-1,\ell)$ or 
$(r,s)=(\ell-1,\ell-1)$.
By Poincar\'e duality, there are elements $a,b\in H^*(X)\subset H^*(M)$ so 
that $a\cup u = b \cup v=\omega \in H^{2\ell-1}(M)$.  

If $(r,s)=(\ell-1,\ell)$, then $a\in H^\ell(X)$ and $b\in H^{\ell-1}(X)$.  
Then, since $X$ is smooth, $W_{\ell-2}(H^{\ell-1}(X))=0$, and 
$b$ is of weight at 
least $\ell-1$.  Since $W_{\ell+1}(H^\ell(X))=0$ by hypothesis, $a$ is of 
weight at least $\ell+2$.  Since $\omega$ is of weight $2\ell$, is follows 
that $u$ is of weight at most $\ell-2$ and $v$ is of weight at most $\ell+1$.  
Consequently, $u \cup v $ is of weight at most $2\ell-1$ in $H^{2\ell-1}(M)$, 
which is pure of weight $2\ell$.  Hence $u\cup v=0$.  
If $(r,s)=(\ell-1,\ell-1)$, then $a,b\in H^{\ell}(X)$ are both of weight at 
least $\ell+2$, and a similar argument shows that $u \cup v=0$.
\qed\end{pf}

\subsection{Plane algebraic curves}
\label{subsec: plane curves}

For an arbitrary projective hypersurface, the cohomology ring
of the boundary manifold (with $\C$ coefficients) need not 
admit the structure of a double.
We illustrate this phenomenon in dimension two.

\begin{theorem}
\label{thm:bdry curve}
Let $V=V_1\cup \dots\cup V_k$ be a reducible algebraic
curve in $\CP^2$, with complement $X$ and boundary 
manifold $M$.  Then $H^*(M;\C) \cong \sD(H^*(X;\C))$ 
if and only if all the irreducible components $V_j$ are 
rational curves.
\end{theorem}

\begin{pf}
If an irreducible component $V_j$ of $V$ is a rational curve, then the
normalization of $V_j$ has genus $0$.  It follows that all nontrivial
cohomology classes in $H^1_0(V_j)=H^1(V_j)$ are of weight $0$.  Using
this, an inductive argument with the Mayer-Vietoris sequence reveals that
the same holds for $H^1_0(V)=H^1(V)$.  It follows that $H^2(X) \cong
H^1_0(V)$ is pure of weight $4$, see \cite[p.~246]{Dimca}.  So $H^*(M)
\cong \sD(H^*(X))$ by Theorem \ref{thm:weight double}.

Conversely, if an irreducible component $V_j$ of $V$ is not 
a rational curve, then the degree of $V_j$ is necessarily at 
least three.  In this situation, $H^1(V)=H^1_0(V)$ contains 
nontrivial classes of weights $0$ and $1$ (see \cite{Dur83}).
It follows that $H^2(X)$ contains classes of weights $3$ 
and $4$ (see \cite{Dimca}).  (Note that the weight condition 
of Theorem \ref{thm:weight double} fails.)  In this instance, 
it is readily checked that the cup-product 
$H^1_0(V) \otimes H^1_0(V)\to H^2_0(V)$ is nontrivial.  
Hence, $H^*_0(V) \subset H^*(M)$ is not a square-zero 
subalgebra, compare Theorem \ref{thm:square-zero}, 
and $H^*(M) \not\cong \sD(H^*(X))$.
\qed\end{pf}

Suppose $V$ is an arrangement of rational curves in 
$\CP^2$, with complement $X$, and boundary manifold $M$.  
A presentation for the cohomology ring $H^*(X;\C)$ was given 
in \cite[Theorem 0.4]{Cog}.  Our Theorem~\ref{thm:bdry curve} 
can now be used to compute the cohomology ring $H^*(M;\C)$.  

\section{Hyperplane arrangements} 
\label{sec:arrangements}

Let $\A$ be an arrangement of hyperplanes in $\CP^{\ell}$.  
For each hyperplane $H$ of $\A$, let $f_H$ be a linear form 
with $H=\{f_H=0\}$.  Then $f=Q(\A)=\prod_{H\in\A} f_H$ is a 
defining polynomial for $\A$, the
hypersurface $V=V(\A)$ is given by $V = f^{-1}(0) =
\bigcup_{H\in \A} H$, and the complement of the arrangement is
$X=X(\A)=\CP^\ell \setminus V$. 

\subsection{Boundary manifold of an arrangement} 
\label{subsec: bdry arr} 
Let $M=M(\A)$ be the boundary manifold of the hypersurface 
$V=V(\A)$.  The next theorem expresses the (integral) cohomology 
ring of $M$ in terms of the Orlik-Solomon algebra 
$A=A(\A)=H^*(X(\A);\Z)$ of the arrangement $\A$.  

\begin{theorem}
\label{thm:coho bdry arr}
Let $\A$ be an arrangement of hyperplanes in $\CP^\ell$ 
with complement $X$ and boundary manifold $M$.  
Then $H^*(M;\Z) \cong \sD(H^*(X;\Z))$.
\end{theorem}

\begin{pf}
For any hyperplane arrangement $\A$, the cohomology $H^k(X,\C)$ 
is pure of weight $2k$, that is, the weight filtration takes the form 
$0=W_{2k-1} \subset W_{2k}= H^k(X;\C)$, for every $k$, 
see Shapiro \cite{Sh93}, and also Kim \cite{Kim}.  
Hence, by Theorem \ref{thm:weight double}, we have 
$H^*(M;\C) \cong \sD(H^*(X;\C))$.

Let $A=H^*(X;\Z)$ be the integral Orlik-Solomon algebra of $\A$.  It is 
well known that $A=\bigoplus_{k=0}^\ell A^k$ is torsion-free.  Let 
$\sD(A)=A\ltimes \bar{A}$ be the integral double of $A$, the 
trivial extension of $A$ by 
$\bar{A}=\bigoplus_{k=\ell-1}^{2\ell-1} \Hom_\Z(A^{2\ell-k-1},\Z)$, 
with $A$-bimodule structure as given in \S\ref{subsec:double}.  Since 
$A=H^*(X;\Z)$ is torsion-free, $H^*(M;\Z)$ is also torsion-free, see 
Proposition \ref{prop:cohomology groups}.  Since 
$H^*(M;\C) \cong \sD(H^*(X;\C))$, it follows that 
$H^*(M;\Z) \cong \sD(A)$.
\qed\end{pf}

Let $L(\A)$ be the intersection poset of the arrangement $\A$, 
the set of all nonempty intersections of elements of $\A$, 
ordered by reverse inclusion.  By the Orlik-Solomon theorem 
(see \cite{OT1,Yuz}), the integral cohomology ring of $X(\A)$ 
is determined by $L(\A)$.  Our next result shows that the 
cohomology of $M(\A)$ is determined by $L(\A)$ and the 
ambient dimension.  

\begin{corollary} 
\label{cor:coho bdry comb}
If $\A_1$ and $\A_2$ are hyperplane arrangements in 
$\CP^{\ell}$ with isomorphic intersection posets,
then $H^*(M(\A_1);\Z)\cong H^*(M(\A_2);\Z)$.
\end{corollary}
\begin{pf}
By the Orlik-Solomon theorem, 
if $L(\A_1)\cong L(\A_2)$, then 
$A(\A_1) \cong A(\A_2)$. 
Proposition \ref{prop:functorial} implies that the (integral) doubles 
are isomorphic.  Thus, by Theorem \ref{thm:coho bdry arr}, 
$H^*(M(\A_1);\Z)\cong H^*(M(\A_2);\Z)$.
\qed\end{pf}

\subsection{Computing cup products}
\label{subsec:cup products}
We now exhibit an explicit basis for the cohomology of the boundary
manifold of an arrangement, and compute cup products in that basis.  
Write $\A=\{H_0,H_1,\dots,H_n\}$, and 
designate $H_0$ as the hyperplane at infinity in $\CP^\ell$.  
Let $\A' = \{H_1,\dots,H_n\}$ be the corresponding affine 
arrangement in $\C^\ell = \CP^\ell \setminus H_0$.  
Notice that $\A$ is the projective closure of $\A'$.   

The rank of the affine arrangement $\A'$ is the maximal number 
of linearly independent hyperplanes in $\A'$.  If $\A' \subset \C^\ell$ 
has rank $\ell$, then $\A'$ is said to be essential. 
Observe that the projective arrangement $\A \subset \CP^\ell$ is 
essential if it contains $\ell+1$ independent hyperplanes.  
For an arrangement of rank $r$, it is well known that 
the Betti numbers, $b_k(X)$, of the complement are 
nonzero for all $k$, $0\le k \le r$.  See \cite{OT1} as 
a general reference.

Order the hyperplanes of
$\A'=\{H_1,\dots,H_n\}$ by their indices.  A circuit is an
inclusion-minimal dependent set of hyperplanes (in $\A'$), and a broken
circuit is a set $S$ for which there exists $j < \min(S)$ so that
$\{H_j\} \cup\{H_i \mid i \in S\}$ is a circuit.  Let $\nbc=\nbc(\A')$
denote the collection of subsets $I \subset [n]$ for which $\bigcap_{i
\in I} H_i \neq \emptyset$ and $I$ contains no broken circuits.  
If the rank of $\A'$ is $r$, then all elements of $\nbc$ are of cardinality 
at most $r$.  Note also that $\emptyset \in \nbc$.

 Clearly, the complement of $\A$ in $\CP^\ell$ is diffeomorphic to the 
complement of  $\A'$ in $\C^\ell$.  
The integral cohomology of $X=X(\A)=X(\A')$ is isomorphic to the
Orlik-Solomon algebra $A=A(\A')$, a quotient of an exterior 
algebra on $n$ generators in degree $1$.  A basis for $A$ 
is indexed by the set $\nbc$; denote
this basis for $A$ by $\{a_I \mid I \in \nbc\}$.  If $|I|=k$, then
$a_I \in A^k$.  In particular, the unit in $A$ is $1=a_\emptyset \in
A^0$.  Express the cup-product in $A=H^*(X)$ by
\begin{equation} \label{eq:cupX}
a_I a_J = \sum_{K\in \nbc} \mu_{I,J,K}a_K.
\end{equation}

Denote the images of the generators $a_I$ of $A=H^*(X)$ in 
$H^*(M)$ by the same symbols.  By Poincar\'e duality, there are 
elements $\bar a_I \in H^*(M)$ so that $a_I \bar a_J = \delta_{I,J}\omega$, 
where $\omega$ is a (fixed) generator of $H^{2\ell-1}(M)\cong \Z$.  
In particular, $\bar a_\emptyset=\omega$.  Since $H^*(M)=\sD(A)$, using 
\eqref{eq:multiplication}, 
we obtain the following.

\begin{corollary} 
\label{cor:arr cohomology ring}
The set $\{a_I, \bar a_I \mid I \in \nbc\}$ forms a basis for
$H^*(M)$, and the cup-product in $H^*(M)$
is given by
\begin{equation*} 
\label{eq:cupM}
a_I a_J = \sum_{K\in \nbc} \mu_{I,J,K}a_K, \quad
a_J \bar a_K = \sum_{I\in \nbc} \mu_{I,J,K}\bar a_I, 
\quad \bar a_I \bar a_J = 0.
\end{equation*}
\end{corollary}

\begin{example} 
\label{ex:pencil}
{\rm
Let $\A$ be a near-pencil
in $\CP^2$, defined by the polynomial $Q(\A)=x_0(x_1^n-x_2^n)$.  
As noted in Example \ref{ex:boundary near pencil}, the 
boundary manifold $M$ is diffeomorphic to $\Sigma_{n-1} \times S^1$.

The complement $X$ of $\A$ has Poincar\'e polynomial 
$P(X,t)=1+nt+(n-1)t^2$.  
The \textbf{nbc} basis of the Orlik-Solomon algebra 
$A=H^*(X)$ is given by 
$\{1=a_\emptyset,a_1,\dots,a_n,a_{1,2},\dots,a_{1,n}\}$.  
The cup-product in $A$ is given by $a_1 a_j=a_{1,j}$ and $a_i a_j=
a_{1,j} - a_{1,i}$ for $i>1$. 

The boundary manifold $M$ has Poincar\'e polynomial
$P(M,t)=1+(2n-1)t+(2n-1)t^2+t^3$.  A basis for the cohomology 
ring $\sD(A)=H^*(M)$ is given by the above basis
for the Orlik-Solomon algebra, together with the dual classes
$\{\bar a_{1,2},\dots,\bar a_{1,n},\bar a_1,\dots,\bar a_n,\bar 
a_\emptyset=\omega\}$.  By Corollary \ref{cor:arr cohomology ring}, 
the cup-product in $\sD(A)$ is given by the multiplication in $A$ 
recorded above, $\bar a_I \bar a_J =0$ for all $I$ and $J$, 
$a_j \bar a_k =  a_{1,j} \bar a_{1,k} = \delta_{j,k}\omega$, and 
$a_j \bar a_{1,k} = -\bar a_k+\delta_{j,k}(\bar a_1+\dots+\bar a_n)$.

Now, $H^*(\Sigma_{n-1} \times S^1)=
H^*(\Sigma_{n-1}) \otimes H^*(S^1)$ is generated by $\alpha_j \otimes
1$, $\beta_j \otimes 1$, $1\le j\le n-1$, $\Gamma \otimes 1$, and
$1\otimes z$, where $\alpha_j$, $\beta_j$, $\Gamma$ generate
$H^*(\Sigma_{n-1})$ and satisfy $\alpha_j \beta_k = \delta_{j,k}\Gamma$, and
$z$ generates $H^*(S^1)$.   An explicit isomorphism from 
$H^*(\Sigma_{n-1} \times S^1)$ to $\sD(A)$ is defined by 
\[ 
\alpha_j\otimes 1 \mapsto
a_{j+1}-a_1, \ 
\beta_j \otimes 1 \mapsto \bar a_{1,j+1}, \ 
1 \otimes z \mapsto a_1, \ 
\Gamma \otimes 1 \mapsto \bar a_1+\dots+\bar a_n.  
\] 
}
\end{example}  

\subsection{The $K(\pi,1)$ problem} 
\label{subsec:discuss arr} 

A hyperplane arrangement $\A$ is said to be a $K(\pi,1)$-arrangement 
if the complement $X=X(\A)$ is aspherical, i.e., its universal cover  
is contractible.   Classical examples include the  braid arrangement 
(Fadell-Neuwirth), certain reflection arrangements (Brieskorn) and 
simplicial arrangements (Deligne).  

The boundary manifold of an arrangement in $\CP^1$ is a disjoint 
union of circles.  For $\ell \ge 3$, 
Proposition \ref{prop:nonaspherical} shows that the boundary manifold 
of an arrangement in $\CP^\ell$ is never aspherical.  In the remaining 
case, $\ell=2$, we have the following result.

\begin{proposition} 
\label{prop:aspherical} 
Let $\A$ be a line arrangement in $\CP^2$.  The boundary manifold 
$M=M(\A)$ is aspherical if and only if $\A$ is essential.  
\end{proposition}

\begin{pf}
If $\A$ is not essential, then $\A$ is a pencil of lines in $\CP^2$, 
and so, by Example \ref{ex:boundary pencil}, $M$ is a connected 
sum of $S^1\times S^2$'s.  Thus, $\pi_2(M)\ne 0$. 

If $\A$ is essential, it follows from work of Jiang and Yau \cite{JY98}
that $M$ is an irreducible, sufficiently large Waldhausen graph 
manifold.  Hence, $M$ is aspherical. (In fact,  by \cite{Sch}, 
$M$ admits a metric of non-positive curvature.)
\qed\end{pf}

\section{Topological complexity}
\label{sec: top complexity} 

In this section, we relate the topological complexity of the boundary 
manifold of a hyperplane arrangement to that of the complement. 
We start by relating the zero-divisor length of a graded algebra 
to that of its double.  

\subsection{Cup length and zero-divisor length}
\label{subsec: length} 
Let $A=\bigoplus_{k=0}^{\ell}A^k$  be a graded algebra 
over a field $\k$  (as usual, we assume all graded pieces 
are finite-dimensional).  Define the {\em cup length} of $A$, 
denoted $\cl(A)$, to be the largest integer $q$ for which there 
exist homogeneous elements $a_1, \dots ,a_q\in A^{>0}$ such that 
$a_1\cdots a_q\ne 0$. 

The tensor product $A\otimes A$ has a natural  graded algebra 
structure, with multiplication given by $(u_1\otimes v_1)\cdot 
(u_2\otimes v_2) = (-1)^{|v_1|\cdot|u_2|}u_1u_2\otimes v_1v_2$.
Multiplication in $A$ defines an algebra homomorphism
$\mu\colon A\otimes A\to A$.  Let $J(A)$ be the kernel of this map.
The {\em zero-divisor length} of $A$, denoted by $\zcl(A)$, 
is the length of the longest non-trivial product in this ideal.

\begin{lemma}
\label{lem:J generators}
The ideal $J(A)=\ker(\mu\colon A\otimes A \to A)$ is generated by 
the set $\{\zeta_a:=a\otimes 1 - 1\otimes a \mid a \in A\}$.
\end{lemma}
\begin{pf}
Let $z=\sum_{i=1}^k a_i \otimes b_i$ be an element of $J(A)$.  
Then $\sum_{i=1}^k a_i b_i=0$ in $A$, and it is readily checked 
that $z-\sum_{i=1}^k \zeta_{a_i}(1\otimes b_i)=
\sum_{i=1}^k 1\otimes a_i b_i=1\otimes(\sum_{i=1}^k a_i b_i)=0$ 
in $A\otimes A$.
\qed\end{pf}

These two notions of length behave quite nicely with respect to 
the doubling operation for graded algebras.

\begin{proposition}
\label{prop:length}
Let $A$ be a connected, finite-dimensional graded algebra, with double 
$\sD(A)=A\ltimes \bar{A}$.  Then, $\cl(\sD(A))=\cl(A)+1$ and 
$\zcl(\sD(A))= \zcl(A)+2$.
\end{proposition}

\begin{pf}
Suppose that $\cl(A)=q$, and let $a=a_1\cdots a_q$ be an element in $A$ of
length $q$.  Then $a\cdot \bar{a}$ is a nonzero element in $\sD(A)$, of
length $q+1$.  Thus $\cl(\sD(A)) \ge \cl(A)+1$.  The equality
$\cl(\sD(A))=\cl(A)+1$ then follows from the fact that
$\bar{A}$ is a square-zero ideal in $\sD(A)$.

Next, suppose that $\zcl(A)=q$, and let $z=z_1\cdots z_q$ be an element in  
$J(A)$ of length $q$.  Recall the basis $\{a^k_j\}$ of $A$ from 
\S\ref{subsec:double}, and write
$z= \sum c^{k_1,k_2}_{j_1,j_2} a^{k_1}_{j_1} \otimes a^{k_2}_{j_2}$.  Let $m$ 
be maximal so that $i_1+i_2=m$ and there is a nonzero coefficient 
$c^{i_1,i_2}_{r_1,r_2}$ in this sum.
Then, one can check that
\[
z(\bar{a}^{i_1}_{r_1}\otimes 1)(1\otimes \bar{a}^{i_2}_{r_2}) = 
\pm c^{i_1,i_2}_{r_1,r_2}\omega\otimes\omega + z',
\]
where $\omega=\bar{1}$ generates $\sD(A)^{2\ell-1}$ and $z'$ is a linear 
combination of elements 
$a^{k_1}_{j_1}\bar{a}^{i_1}_{r_1}\otimes a^{k_2}_{j_2}\bar{a}^{i_2}_{r_2}$ in 
$\sD(A)$ of bidegree 
different from $(2\ell-1,2\ell-1)$.  So $\hat{z}=z(\bar{a}^{i_1}_{r_1}\otimes 
1)(1\otimes \bar{a}^{i_2}_{r_2})$ is a nonzero element in $J(\sD(A))$, 
of length at least $q+2$.  Thus $\zcl(\sD(A)) \ge \zcl(A)+2$.

To show that $\zcl(\sD(A)) = \zcl(A)+2$, it suffices to check that $\hat{z}
\zeta_\alpha= \hat{z}(\alpha\otimes 1-1\otimes\alpha)=0$ for  
$\alpha\in\sD(A)$.  We may assume that $\alpha$ is an element of the basis
$\{a^k_j,\bar{a}^k_j\}$ for $\sD(A)$.  If $\alpha=\bar{a}^k_j\in\bar{A}$, then 
$\hat{z} \zeta_\alpha=0$ since $\bar{A}$ is a square-zero ideal in $\sD(A)$.  
If $\alpha=a^k_j\in A$ and $\hat{z} \zeta_\alpha \neq 0$, then 
$z \zeta_\alpha$ is a nonzero element of length $q+1$ in $J(A)$, 
contradicting the assumption that $\zcl(A)=q$.
\qed\end{pf} 

\subsection{LS category and topological complexity}
\label{subsec: ls and tc} 

Let $p\colon Y\to X$ be a fibration.  The 
{\em sectional category} of $p$, denoted $\secat(p)$, 
is the smallest integer $q$ such that $X$ can be covered 
by $q$ open subsets, over each of which $p$ has a section.  
A cohomological lower bound is given by:  
\begin{equation}
\label{eq:secat}
\secat (p)> \cl (\ker (p^*\colon H^*(Y;\k)\to H^*(X;\k))),
\end{equation}
see James \cite{Ja} as a classical reference.
If $p\colon PX\to X$ is the path-fibration of a pointed space $X$, 
then $\secat(p)=\cat(X)$, the Lusternik-Schnirelmann category 
of $X$.  The category of $X$ depends only on the homotopy-type 
of $X$.  Since $PX$ is contractible, the inequality \eqref{eq:secat}
reduces to  $\cat(X) > \cl(X):=\cl(H^*(X;\k))$. If $X$ is a finite 
simplicial complex, then $\cat(X)\le \dim (X)+1$. Furthermore, 
$\cat(X\times Y)\le \cat(X) + \cat(Y)-1$.  

Now let $X^I$ be the space of all continuous paths from 
$I=[0,1]$ to $X$, with the compact-open topology, 
and let $\pi\colon X^I \to X \times X$ be the fibration 
given by $\pi(\gamma)=(\gamma(0),\gamma(1))$. 
The \emph{topological complexity} of $X$, introduced 
by Farber in \cite{Fa03} and denoted by $\tc(X)$, may 
be realized as the sectional category of $\pi$. Again, 
$\tc(X)=\secat(\pi)$ depends only on the homotopy type of $X$.  
Using the fact that $X^I\simeq X$, and the K\"unneth formula, 
\eqref{eq:secat} reduces to 
$\tc(X)> \zcl(X):=\zcl(H^*(X;\k))$.  If $X$ is a finite 
simplicial complex, then $\cat(X)\le \tc(X) \le 2\cat(X)-1$; 
in particular, $\tc(X)\le 2\dim (X) +1$. 
Furthermore, $\tc(X\times Y)\le \tc(X) + \tc(Y)-1$.  

As noted in \cite{Fa03}, topological complexity 
is {\em not} determined by the LS category.  
For example, $\cat(S^n)=2$ for all $n\ge 1$, whereas $\tc(S^n)=2$ 
for $n$ odd and $\tc(S^n)=3$ for $n$ even; also, $\cat(T^n)=\tc(T^n)=n+1$, 
but $\cat(\Sigma_g)=3$ and $\tc(\Sigma_g)=5$ for $g\ge 2$.

In \cite{FY}, Farber and Yuzvinsky study the invariants $\tc(X)$ and 
$\zcl(X)$ in the case when $X$ is the complement of a 
(central, essential) hyperplane arrangement in $\C^{\ell}$.  
They show that $\tc(X)\le 2\ell$, and that this upper bound 
is attained for some classes of arrangements, including 
generic arrangements of sufficiently large cardinality 
and the reflection arrangements of types A, B, and D.

\subsection{Topological complexity of the boundary manifold}
\label{subsec: tc bdry} 

Using Theorem \ref{thm:coho bdry arr} and 
Proposition \ref{prop:length}, we see that  
the cup and zero-divisor lengths of the boundary 
manifold of an arrangement are determined in a 
simple fashion by the respective lengths of the complement. 

\begin{corollary}
\label{cor:length M and X}
Let $\A$ be an arrangement of hyperplanes in 
$\CP^{\ell}$, with complement $X$ and boundary manifold $M$. 
Then:
\begin{equation*}
\label{eq:length M and X}
\cl(M)=\cl(X)+1 \ \text{and}\  \zcl(M)=\zcl(X)+2.
\end{equation*}
Moreover, if $\A$ is essential, then $\cl(M)=\ell+1$.
\end{corollary}

The relationship between the LS category and topological 
complexity of the boundary manifold on one hand, and the 
complement on the other hand, is more subtle, as the following  
example indicates.  

\begin{example}
\label{ex:boolean tc}
{\rm
Let $\A$ be the Boolean arrangement $\CP^{\ell}$.  
Then $X\simeq T^{\ell}$ and $M=T^{\ell}\times S^{\ell-1}$.  
An easy computation shows that $\cat(M)=\cat(X)+1=\ell+2$; 
on the other hand, $\tc(M)=\tc(X)+2=\ell+3$ if $\ell$ is even, 
but $\tc(M)=\tc(X)+3=\ell+4$ if $\ell$ is odd.
}
\end{example}

For projective line arrangements, we can narrow down 
the possible values of the category and topological complexity 
of the boundary manifold.

\begin{proposition} 
\label{prop:cat tc 3} 
Let $\A$ be a line arrangement in $\CP^2$, with 
boundary manifold $M$.  If $\A$ is not essential, then 
$\cat(M)=2$ or $3$ and  $\tc(M)=4$, $5$, or $6$.  
If $\A$ is essential, then 
$\cat(M)=4$ and $\tc(M)=5$, $6$, or $7$.
\end{proposition}

\begin{pf}
As shown in \cite{GG}, the LS category of a closed $3$-manifold 
$M$ depends only on $\pi_1(M)$: it is $2$, $3$, or $4$, according 
to whether $\pi_1(M)$ is trivial, a non-trivial free group, or not a 
free group.  

Suppose $\A$ is a pencil of $n+1$ lines.  If $n=0$, then $M=S^3$, 
so $\cat(M)=\tc(M)=2$.  If  $n=1$, then $M=S^1\times S^2$, 
so $\cat(M)=3$ and $\tc(M)=4$.  If $n>1$, then 
$M=\#^n S^1\times S^2$, and so $\cat(M)=3$, and $\tc(M)= 5$.

On the other hand, if $\A$ is essential, then, as noted in Proposition 
\ref{prop:aspherical}, $M$ is aspherical. In particular, $\cd(\pi_1(M))=3$, 
and so $\pi_1(M)$ cannot be free.  Hence, $\cat(M)=4$, and the bounds 
on $\tc(M)$ follow at once.
\qed\end{pf}

All the various possibilities listed in Proposition \ref{prop:cat tc 3} 
do occur.  For example, if $\A$ is a near-pencil of $n+1\ge 4$ 
lines, then $M=\Sigma_{n-1} \times S^1$, and so 
$\cat(M)=4$ and $\tc(M)=6$.   
We summarize in Table \ref{table:cat and tc} the possible values 
for the LS category and topological complexity of both the 
complement and the boundary manifold of an arrangement in 
$\CP^2$, together with sample representatives for the defining 
polynomials.  

\renewcommand{\arraystretch}{1.4}
\begin{table}
\[
\begin{array}{|l|c|c|c|c|c|c|}
\hline
\cat(X) & 1 & 2& 2 & 3 & 3 & 3\\
\hline
\tc(X) & 2 & 3 &4 & 4 & 5 & 6\\
\hline
\cat(M) & 2 & 3& 3 & 4 & 4 & 4\\
\hline
\tc(M) & 2 & 4 &5 & 5 & 6 & 7\\
\hline
f & x_0 & x_0x_1 & x_0^3-x_1^3 
& x_0x_1x_2 & x_0(x_1^3-x_2^3) 
&  (x_0^2-x_1^2)(x_0^2-x_2^2)(x_1^2-x_2^2) \\
\hline
\end{array}
\]
\caption{Possible values of LS category and topological 
complexity for the complement $X$ and boundary manifold $M$ 
of a line arrangement in $\CP^2$.}
\label{table:cat and tc}
\end{table}
\renewcommand{\arraystretch}{1.0}

In high dimensions, a complete understanding of the possible 
values for $\cat(M)$ and $\tc(M)$ is not at hand.  
Nevertheless, we have the following class of arrangements 
(mentioned in Example~\ref{ex:boundary product pencils}), 
where precise formulas can be given.  

\begin{proposition}
\label{ex:product pencil tc}
Let $\A$ be the hyperplane arrangement in $\CP^{2k}$ defined 
by the polynomial 
$f=x_0  \prod_{i=1}^{k} (x_i^{n_i} - y_i^{n_i})$, 
with  $n_i\ge 3$.  
If $X$ is the complement and $M$ is the boundary manifold, then:
\begin{align*}
\cat(X)&=2k+1, 
& \tc(X)&=3k+ 1,  \\
\cat(M)&=2k+2, 
& \tc(M)&=3k+3. 
\end{align*}
\end{proposition}

\begin{pf}
We have $X\simeq T^k \times  \prod_{i=1}^k \bigvee^{n_i-1} S^1$, 
while $M=T^k \times (\#^{m} T^k \times S^{2k-1})$, where 
$m=\prod_{i=1}^{k} (n_i-1)$.  A computation 
shows that $\cl(X)=2k$ and $\zcl(X)=3k$.  
Hence, by Corollary \ref{cor:length M and X}, 
$\cl(M)=2k+1$  
and $\zcl(M)=3k+2$.    

Let $W=T^k \times S^{2k-1}$.  Note that 
$\cl(W)=k+1$, while $\cat(W)\le \cat(T^k)+\cat (S^{2k-1})-1=k+2$; hence 
$\cat(W)=k+2$.  In fact, if we consider $W$ with its standard 
CW decomposition, we can take $W=\bigcup_{i=0}^{k+1} U_i$, 
with $U_0$ a small ball around the $0$-cell $e^0$, $U_i$ the 
union of the (open) $i$ and $i+2k-1$ cells, for $1\le i\le k$, 
and $U_{k+1}$ the top cell $e^{3k-1}$; plainly, each $U_i$ is 
contractible in $W$. 
 
Now, $\#^m W$ is obtained by attaching a top cell to the wedge 
of $m$ copies of $W\setminus U_{k+1}$ at the basepoint $e^0$; 
thus, we may find a decomposition $\#^m W=\bigcup_{i=0}^{k+1} V_i$ 
as before, with $V_i$ contractible in $\#^m W$.  It follows that 
$\cat(\#^m W)=k+2$, and so $\tc(\#^m W)\le 2k+3$.  Thus, 
$\tc(M)\le \tc(T^k)+ \tc(\#^m W)-1 \le 3k+3$, and we are done.  
\qed\end{pf}

As a consequence, we see that the difference between 
the topological complexity and the LS category of the boundary 
manifold of an arrangement can be arbitrarily large.

\begin{corollary}
\label{cor:tc-cat}
For each $k\ge 1$, there is an arrangement $\A$ with boundary 
manifold $M=M(\A)$ for which $\tc(M)-\cat(M)=k$. 
\end{corollary}

\section{Resonance}
\label{sec: res boundary}

In this section, we study the resonance varieties of the 
trivial extension of a graded algebra.  As an application, 
we obtain information about the structure of the resonance 
varieties of the boundary manifold of a hyperplane arrangement.  
Throughout, let $\k$ be an algebraically closed field of 
characteristic~$0$.

\subsection{Resonance varieties}
\label{subsec: res var}

Let $A=\bigoplus_{k=0}^{\ell} A^k$ be a graded, 
graded-commutative, connected algebra over $\k$.   
Assume each graded piece $A^k$ is finite-dimensional.  
For each $a\in A^1$, we have $a\cdot a=0$; thus, 
multiplication by $a$ defines a cochain complex
\begin{equation}
\label{eq:aomoto}
\xymatrixcolsep{22pt}
(A,a)\colon \:\:
\xymatrix{
0 \ar[r] &A^0 \ar[r]^{a} & A^1
\ar[r]^{a}  & A^2 \ar[r]^(.46){a}&\, \cdots \, \ar[r]^{a}
& A^{\ell}\ar[r] & 0}.
\end{equation}

By definition, the {\em resonance varieties} of $A$ are the 
jumping loci for the cohomology of these complexes:
\begin{equation} 
\label{eq:res var}
\RR^{k}_{d}(A)=\{ a\in A^1 \mid \dim_{\k} H^k(A,a) \ge d\},
\end{equation}
for $0\le k\le \ell$ and $0\le d \le b_k=b_k(A)$.  Notice that 
$A^1=R^k_0(A)\supset R^k_1(A)\supset \cdots \supset 
R^k_{b_k}(A) \supset \{0\}$. 
The sets $\RR^{k}_{d}(A)$ are algebraic subvarieties of the affine
space $A^1 = \k^n$, and are isomorphism-type invariants 
of the graded algebra $A$. They have 
been the subject of considerable recent interest, particularly 
in the context of hyperplane arrangements, see for instance 
\cite{CScv,MS,Yuz}, and references therein.

An element $a\in A^1$ is said to be {\em nonresonant} 
if the dimensions of the cohomology groups $H^*(A,a)$ 
are minimal.  If $A$ is the Orlik-Solomon 
algebra of an arrangement of rank $\ell$, and $a\in A^1$ 
is nonresonant, then $H^k(A,a)=0$ for $k\ne \ell$, 
see for instance \cite{Yuz}.

\subsection{Resonance varieties of a doubled algebra}
\label{subsec: res var double}

We now compare the resonance varieties of $A$ to those 
of the doubled algebra $\sD(A)=A\ltimes \bar{A}$, under the 
assumption that $\ell\ge 3$. Notice that for such $\ell$, we 
have $\sD(A)^1 = A^1$.

\begin{theorem} 
\label{thm:boring res}
If $A$ is a graded, connected $\k$-algebra  
and $\ell \ge 3$, then the resonance
varieties of $\sD(A)$ are given by
\begin{equation*}
\RR^{k}_{d}(\sD(A))=
\begin{cases}
\RR^{k}_{d}(A)
& \text{if $k \le \ell-2$,} \\[3pt]
\bigcup_{p+q=d} \bigl(\RR^{\ell-1}_{p}(A) \cap
\RR^{\ell}_{q}(A)\bigr)
& \text{if $k=\ell-1$ or $k=\ell$,} \\[3pt]
\RR^{2\ell-k-1}_{d}(A)
& \text{if $k \ge \ell+1$.}
\end{cases}
\end{equation*}
\end{theorem}

\begin{pf}  
Fix a basis $\{a^p_i\}$ for $A$, and let $\{a_i^p,\bar a_j^q\}$ 
be the corresponding basis for the double $\sD(A)$ 
as in \S\ref{subsec:double}.
Let $m_k=m_k(a)$ and $\sD(m_k)=\sD(m_k(a))$ denote the matrices
of the maps $A^{k-1} \to A^k$ and $\sD(A)^{k-1} \to \sD(A)^k$
given by multiplication by $a \in A^1 = \sD(A)^1$ in the bases specified 
above.  
An exercise
in linear algebra reveals that
\begin{equation*}
\RR^{k}_{d}(A)=\set{a\in A^1 \mid
\dim_{\k} A^k - \rank m_k - \rank m_{k+1}\ge d}.
\end{equation*}
Similarly,
\begin{equation*}
\RR^{k}_{d}(\sD(A))=\set{a\in \sD(A)^1 \mid
\dim_{\k} \sD(A)^k - \rank \sD(m_k) - \rank \sD(m_{k+1})\ge d}.
\end{equation*}

The complex $(\sD(A),a)$ has terms
$\sD(A)^k=A^k$ for $k\le \ell-2$,
$\sD(A)^{\ell-1}=A^{\ell-1} \oplus \bar{A}^\ell$,
$\sD(A)^\ell=\bar{A}^{\ell-1} \oplus A^\ell$, and
$\sD(A)^k=\bar{A}^{2\ell-k-1}$ for $k \ge \ell+1$.
Using the definition of the multiplication in $\sD(A)$, one
can check that the boundary maps of this complex have
matrices
$\sD(m_k)=m_k$ for $k\le\ell-2$,
\begin{equation*}
\sD(m_{\ell-1})=\left[\begin{matrix}m_{\ell-1} &
0\end{matrix}\right],
\quad
\sD(m_\ell)=
\left[\begin{matrix}0 & m_\ell \\ \pm m_\ell^\tp &
0\end{matrix}\right],
\quad
\sD(m_{\ell+1})=
\left[\begin{matrix}\pm m_{\ell-1}^\tp \\  0\end{matrix}\right],
\end{equation*}
and $\sD(m_k)=\pm m_{2\ell-k}$ for $k \ge \ell+2$.
Calculating ranks of these matrices, and using the above
descriptions of the resonance varieties $\RR^{k}_{d}(A)$ 
and $\RR^{k}_{d}(\sD(A))$ yields the result.
\qed\end{pf}

If $\ell=2$, then $\sD(A)^1 = A^1 \oplus \bar{A}^2$.  
If $(a,b)\cdot (a,b)=0$ for all $(a,b)\in\sD(A)^1$, then 
$\bigl(\sD(A),(a,b)\bigr)$ is a cochain complex for each 
$(a,b)$ as in \eqref{eq:aomoto}, and the resonance 
varieties of $\sD(A)$ are
\[
\RR_d^k(\sD(A))=\{(a,b)\in \sD(A)^1 \mid 
\dim_\k H^k(\sD(A),(a,b)) \ge d\}.
\]
In this situation, the boundary maps of the complex 
$(\sD(A),(a,b))$ have matrices 
\[
\sD(m_1)=\left[\begin{matrix}m_1&\bar{m}_1\end{matrix}\right],
\quad
\sD(m_2)=\left[\begin{matrix}\phi& m_2\\ -m_2^\tp & 0\end{matrix}\right],
\quad
\sD(m_3)=\left[\begin{matrix}m_1^\tp \\ \bar{m}_1^\tp\end{matrix}\right],
\]
where, as above, $m_k=m_k(a)$ is the matrix of multiplication 
by $a$, $A^{k-1}\to A^k$, $\bar{m}_1=\bar{m}_1(b)$ is the matrix 
of multiplication by $b$, $\bar{A}^2 \to \bar{A}^1$, and $\phi=\phi(b)$ 
is the matrix of multiplication by $b$, $A^1 \to \bar{A}^1$.  Since 
$A$ and $\sD(A)$ are graded commutative, the matrix $\phi$ is 
skew-symmetric.  The structure of these matrices, 
$\sD(m_3)=\sD(m_1)^\tp$ and $\sD(m_2)^\tp = -\sD(m_2)$, follows 
from the multiplication in $\sD(A)$, see \eqref{eq:multiplication}.

\subsection{Aomoto complexes}
\label{subsec:aomoto}
The complex \eqref{eq:aomoto} may be realized as the 
specialization at $a$ of the Aomoto complex of the algebra $A$.  
Let $R_A=\Sym(A_1)$ be the symmetric algebra on the 
$\k$-dual of $A^1$, and let $\mathbf{x}=\{x_1,\dots ,x_n\}$ 
be the basis for $A_1$ dual to the basis $\{a_1^1,\dots,a_n^1\}$ 
for $A^1$.  Then $R_A$ becomes identified with the polynomial 
ring $R=\k[\mathbf{x}]$.  The Aomoto complex of $A$ is the cochain 
complex 
\begin{equation}
\label{eq:aomoto again}
\xymatrixcolsep{22pt}
\xymatrix{
A^0\otimes_\k R \ar[r]^{d^1} & A^1\otimes_\k R
\ar[r]^{d^2}  & A^2\otimes_\k R \ar[r]^(.58){d^3}
&\, \cdots \, \ar[r]^(.44){d^\ell} & A^{\ell}\otimes_\k R},
\end{equation}
where the boundary maps are multiplication by 
$\sum_{j=1}^n a^1_j \otimes x_j$.  Notice that the multiplication map 
$\mu\colon A^1 \otimes A^{p-1} \to A^p$ can be recovered from 
the boundary map $d^p$.  Denote the matrix of $d^1$ by 
$d_{\mathbf{x}}$, and that of $d^2$ by $\Delta=\Delta_A$.  
If the multiplication $A^1\otimes A^1 \to A^2$ is given by 
$a_i^1 a_j^1=\sum_{k=1}^m \mu_{i,j,k} a^2_k$, the latter is 
an $n \times m$ matrix of linear forms over $R$, with entries
\begin{equation}
\label{eq:delta matrix}
\Delta_{j,k}=\sum_{i=1}^n \mu_{i,j,k} x_i.
\end{equation}
The (transpose of the) matrix $\Delta_A$ is the (linearized) 
\emph{Alexander matrix} of the algebra $A$, which appears 
in various guises in, for instance, \cite{CSai}, \cite{CScv}, 
\cite{MS}, \cite{PS}.

The Aomoto complex of the double $\sD(A)$ may be constructed 
analogously.  In light of Theorem \ref{thm:boring res}, we focus 
on the case $\ell=2$.  Here, $\sD(A)^1=A^1 \oplus \bar{A}^2$, 
with basis $\{a^1_i,\bar a^2_j\}$, where $1\le i\le n$ and say 
$1\le j \le m$.  Identify the ring $R_{\sD(A)}=\Sym(A_1 \oplus 
\bar{A}_2)$ with the polynomial ring $\k[\mathbf{x},\mathbf{y}]$.
Then, the Aomoto complex of $\sD(A)$ is the complex  
\begin{equation}
\label{eq:aomoto D}
\xymatrix{
\sD(A)^0\otimes_\k S \ar[r]^{D^1} & \sD(A)^1 \otimes_\k S 
 \ar[r]^{D^2} & \sD(A)^2\otimes_\k S  \ar[r]^{D^3} 
 & \sD(A)^3 \otimes_\k S
},
\end{equation}
where the boundary maps are multiplication by 
$\sum_{i=1}^n a^1_i \otimes x_i + 
\sum_{j=1}^m \bar a^2_j \otimes y_j$.  
Denote the matrix of $D^1$ by 
$\begin{pmatrix}d_{\mathbf{x}}&d_{\mathbf{y}}\end{pmatrix}$, 
and that of $D^2$ by $\Delta_{\sD(A)}$.  Then it follows from 
\eqref{eq:multiplication} that the matrix of $D^3$ is 
$\begin{pmatrix}d_{\mathbf{x}}&d_{\mathbf{y}}\end{pmatrix}^\tp$, 
and that
\begin{equation}
\label{eq:bmat}
\Delta_{\sD(A)}  = \begin{pmatrix} 
\Phi & \Delta_A \\ 
-\Delta_A^\tp & 0 \end{pmatrix}, 
\end{equation}
where $\Phi$ is the $n \times n$ matrix with entries 
$\Phi_{i,j} = \sum_{k=1}^{m} \mu_{i,j,k} y_k$.  
Notice that $\Delta_{\sD(A)}$ is a skew-symmetric matrix 
of linear forms, and that $d_{\mathbf{x}} \Phi=d_{\mathbf{y}}\Delta_A^\tp$.

If $A=\bigoplus_{k=0}^\ell A^k$ and $\ell \ge 3$, the 
relationship between the Aomoto complexes of $A$ 
and $\sD(A)$ is implicit in the proof of Theorem 
\ref{thm:boring res}.  We relate these complexes 
in the case $\ell=2$.  Consider the Aomoto complex 
of $A$ and its dual,
\begin{equation*}
\xymatrixcolsep{22pt}
\xymatrix{
A^0\otimes_\k S \ar[r]^{d_\b{x}} & A^1 \otimes_\k S 
\ar[r]^{\Delta_A} & A^2\otimes_\k S}
\ \ \text{and} \ \  
\xymatrix{
\bar A^2\otimes_\k S \ar[r]^{-\Delta_A^\tp} & \bar A^1 
\otimes_\k S \ar[r]^{-d_\b{x}^\tp} & \bar A^0\otimes_\k S,
}
\end{equation*}
where we have extended scalars and changed the signs 
for reasons which will become apparent.

\begin{lemma} 
\label{lem:mapping cone}
The maps $\{d_\b{y},-\Phi,d_\b{y}^\tp\}$ provide a chain map
\begin{equation*}
\xymatrixcolsep{30pt}
\xymatrix{
A^0\otimes_\k S \ar[r]^{d_\b{x}} \ar[d]^{d_\b{y}}   
& A^1 \otimes_\k S \ar[r]^{\Delta_A} \ar[d]^{-\Phi} 
& A^2\otimes_\k S \ar[d]^{d_\b{y}^\tp} 
\\
\bar A^2\otimes_\k S \ar[r]^{-\Delta_A^\tp} 
& \bar A^1 \otimes_\k S \ar[r]^{-d_\b{x}^\tp}
& \bar A^0\otimes_\k S.
}
\end{equation*}
Furthermore, the Aomoto complex of $\sD(A)$ is the mapping 
cone of this chain map.
\end{lemma}

An alternate way to compute the resonance varieties $\RR^1_{d}(A)$
is by taking the zero locus of the determinantal ideals of 
the linearized Alexander matrix of $A$.    
If $\Psi$ is a $p\times q$ matrix ($p\le q$) with polynomial 
entries, define $\RR_d(\Psi)=V(E_{p-d}(\Psi))$ where 
$E_{r}(\Psi)$ is the ideal of $r\times r$ minors. 
Proceeding as in the proof of Theorem 3.9 from \cite{MS} (see also 
\cite{CScv}), we find that $\RR^1_{d}(A)=\RR_{d}(\Delta_A)$. 
Similarly, $\RR^1_{d}(\sD(A))=\RR_d (\Delta_{\sD(A)})$.  

\subsection{Resonance of arrangements}
\label{subsec:res arr}
For a space $X$ with the homotopy type of a finite CW-complex, 
define the resonance 
varieties of $X$ by $\RR^k_d(X)=\RR^k_d(H^*(X;\k))$.  

Let $\A$ be an arrangement of hyperplanes, with complement 
$X$, boundary manifold $M$, and Orlik-Solomon algebra 
$A=H^*(X;\k)$.  If $\A$ is an arrangement 
in $\CP^\ell$ with $\ell \ge 3$, 
then it follows from Theorem \ref{thm:boring res} that the
resonance varieties of the complement, $\RR^k_d(X)=\RR^k_d(A)$, 
completely determine the resonance varieties, 
$\RR^k_d(M)=\RR^k_d(\sD(A))$, of the boundary manifold.  
So assume that $\A\subset \CP^2$ is a line arrangement.

The complex
$(A,a)$ of \eqref{eq:aomoto} may be realized as the specialization
$\left.A^\bullet \otimes_\k R\right|_{\b{x}\mapsto a}$ of the Aomoto
complex of $A$.  Since $\sD(A)^1 = A^1 \oplus \bar{A}^2$, the resonance
varieties of the boundary manifold are given by 
\[
\RR_d^k(M) = \RR_d^k(\sD(A)) = \{(a,b) \in
A^1\oplus \bar{A}^2 \mid \dim H^k(\sD(A),(a,b)) \ge d\}.
\]  
Similarly, the complex
$(\sD(A),(a,b))$ may be realized as the specialization
$\left.\sD(A)^\bullet \otimes_\k S\right|_{(\b{x},\b{y})\mapsto
(a,b)}$ of the Aomoto complex of $\sD(A)$.  The properties of the
boundary maps of the complex \eqref{eq:aomoto D} noted above imply
that the resonance varieties of $M$ satisfy $\RR_d^k(M) =
\RR_d^{3-k}(M)$.

Recall that, for nonresonant $a\in A^1$, we have $H^k(A,a)=0$ 
for $k=0,1$. Write $b_k=b_k(A)=\dim_\k A^k$, and 
$\beta = 1-b_1 + b_2$. Note that $\beta=\dim_\k H^2(A,a)$.

\begin{proposition} 
\label{prop:res rels}
Let $\A\subset\CP^2$ be a line arrangement with Orlik-Solomon 
algebra $A$ and double $\sD(A)$.
\begin{enumerate}
\item
\label{item:rel1}
If $a \in A^1$ is nonresonant for $A$, then for any $b$, 
$(a,b)\in \sD(A)^1$ is nonresonant 
for $\sD(A)$.  Furthermore, $H^0(\sD(A),(a,b))=H^3(\sD(A),(a,b))=0$
and  $H^1(\sD(A),(a,b))=H^2(\sD(A),(a,b))=\k^\beta$.
\item 
\label{item:rel2}
If $a \in \RR_d^1(A)$ is nonzero, then for any $b$, 
$(a,b) \in \RR_{d+\beta}^1(\sD(A))$.
\item 
\label{item:rel3}
If $b\neq 0$, then $(0,b) \in \RR^1_{d}(\sD(A))$, 
where $d=b_2-1+ \dim_\k \big(\ker \left.\Phi\right|_{\b{y}
\mapsto b}\big)$.
\end{enumerate}
\end{proposition}

\begin{pf}
Given $(a,b) \in \sD(A)^1$, by Lemma~\ref{lem:mapping cone}, there 
is a corresponding short exact sequence of complexes 
$\xymatrixcolsep{16pt}
\xymatrix{0 \ar[r] & (\bar A^\#,a)^{-1}\ar[r] & (\sD(A),(a,b)) \ar[r] & (A,a) 
\ar[r] &0}$:
\begin{equation*}
\xymatrixrowsep{30pt}
\xymatrixcolsep{60pt}
\xymatrix{
& \bar A^2 \ar[d] \ar[r]^{-\Delta_A^{\tp}} 
& \bar A^1 \ar[d] \ar[r]^{-d_\b{x}}
& \bar A^0 \ar[d] 
\\
\sD(A)^0 \ar[d]  \ar[r]^{\mathsmaller{\begin{pmatrix} d_\b{x}&d_\b{y}
\end{pmatrix}}}
& \sD(A)^1 \ar[d] \ar[r]^{\mathsmaller{\begin{pmatrix} 
\Phi & \Delta_A \\ -\Delta_A^\tp & 0\end{pmatrix}}}
&  \sD(A)^2 \ar[d] \ar[r]^{\mathsmaller{\begin{pmatrix}
d_\b{x}^\tp \\ d_\b{y}^\tp\end{pmatrix}}}
& \sD(A)^3 
 \\
A^0 \ar[r]^{d_\b{x}} 
& A^1 \ar[r]^{\Delta_A} 
&  A^2 
}
\end{equation*}
where $(\bar A^\#,a)$ denotes the specialization at $a$ of the dual 
of the Aomoto complex of $A$.  Passing to cohomology yields a 
long exact sequence
\begin{equation}
\label{eq:long exact}
\xymatrixrowsep{5pt}
\xymatrixcolsep{11pt}
\xymatrix{
0 \ar[r] & H^0(\sD(A))\ar[r] & H^0(A) \ar[r] & H^0(\bar A^\#) \ar[r] & 
H^1(\sD(A))
\ar[r] & H^1(A) \ar[r] & \\
\ar[r] &H^1(\bar A^\#)\ar[r] & H^2(\sD(A))
\ar[r] & H^2(A) \ar[r] & H^2(\bar A^\#)\ar[r] & H^3(\sD(A)) \ar[r] & 0
}
\end{equation}
where, for instance $H^k(A)=H^k(A,a)$.  Using the fact that 
$H^k(\bar A^\#)$ is isomorphic to $H^{2-k}(A)$, 
calculations with this long exact 
sequence may be used to establish all three assertions.
\qed\end{pf}
 
As a consequence of Proposition \ref{prop:res rels}, we 
obtain the following.

\begin{corollary}
\label{cor:res var double}
The resonance varieties of the doubled algebra 
$\sD(A)$ satisfy
\begin{enumerate}
\item \label{rr1}
$\RR^1_d(\sD(A)) = \sD^1(A)$ for $d \le \beta$.
\item  \label{rr2}
$\RR^1_{d}(A) \times A^2
\subseteq \RR^1_{d+\beta}(\sD(A))$.
\item  \label{rr3}
$\{0\} \times \RR_{d}(\Phi) 
\subseteq \RR^1_{d+b_2(A)}(\sD(A))$. 
\end{enumerate}
\end{corollary}

All irreducible components of $\RR^1_d(X)=\RR^1_d(A)$ are linear, 
see \cite{CScv}. From items \eqref{rr1} and \eqref{rr2}  in the above 
Corollary, it is clear that the resonance varieties of $M$ contain 
linear components as well. However, item \eqref{rr3} leaves open 
the possibility that $\RR^1_d(M)=\RR^1_{d}(\sD(A))$ contains 
non-linear components, for $d\ge b_2(A)$. This does indeed 
occur, as shown next.

\subsection{General position arrangements}
\label{subsec:gen pos}
Let $\A_n$ be a projective line arrangement consisting of $n+1$ 
lines in general position.  We identify the resonance varieties 
of the boundary manifold $M^3_n=M(\A_n)$. 

The Orlik-Solomon algebra, $A=E/\mathfrak{m}^3$, is the rank 
$2$ truncation of the exterior algebra generated by 
$e_1,\dots, e_n$, where $\mathfrak{m}=(e_1,\dots ,e_n)$.  
Note that $A$ has Betti numbers $b_1=n$,  
$b_2=\binom{n}{2}$, and that $\beta=1-b_1+b_2=\binom{n-1}{2}$. 

For this arrangement, the Alexander matrix $\Delta_A$ is 
the transpose of the matrix of the Koszul differential 
$\delta_2\colon E^2\otimes S \to E^1 \otimes S$.  
The submatrix $\Phi$ of the Alexander matrix $\Delta_{\sD(A)}$ 
recorded in \eqref{eq:bmat} is the generic $n\times n$ skew-symmetric 
matrix of linear forms $\Phi_n$, with entries 
$(\Phi_n)_{i,j}=y_{i,j}$ above the diagonal.

Identity $\sD(A)^1 = E^1 \times E^2$.  Note that 
$\RR^1_d(A)=\RR^1_d(E)=\{0\}$ for $d>0$.  An analysis of the 
long exact sequence \eqref{eq:long exact} in light of this 
observation yields the following sharpening of 
Corollary \ref{cor:res var double} for general 
position arrangements.

\begin{proposition}
\label{prop:res bdry generic}
The resonance varieties of the boundary manifold $M_n$ of 
a general position arrangement of $n+1$ lines in $\CP^2$ 
are given by:
\[
\RR^1_{d}(M_n)=
\begin{cases}
E^1 \times E^2 & \text{if $d\le \binom{n-1}{2}$},\\
\{0\}\times E^2 &\text{if $ \binom{n-1}{2}<d< \binom{n}{2}$},\\
\{0\} \times \RR_{d-\binom{n}{2}} (\Phi_n) &\text{if $\binom{n}{2} \le d <
\binom{n}{2}+n$},\\
\{0\}\times\{0\} &\text{if $d=
\binom{n}{2}+n$}.
\end{cases}
\]
\end{proposition}

If $\Psi$ is a skew-symmetric matrix of size $n$ with 
polynomial entries, define the \emph{Pfaffian varieties} 
of $\Psi$ by
\begin{equation}
\label{eq: pfaff var}
\PP_{d}(\Psi)=V(\Pf_{2\left(\lfloor n/2 \rfloor - d\right)}(\Psi)),
\end{equation}
where $\Pf_{2r}(\Psi)$ is the ideal of $2r \times 2r$ Pfaffians 
of  $\Psi$.  For $n$ even, the ideal $\Pf_n(\Psi)$ is principal, 
generated by $\Pfaff(\Psi)$, the maximal Pfaffian of $\Psi$.  
Well known properties of Pfaffians (see, for instance 
\cite[Cor.~2.6]{BE}) may be used to establish the following 
relationship between the resonance and Pfaffian varieties 
of $\Psi$:
\begin{equation} 
\label{eq:pfaff res}
V(E_{2r-1}(\Psi)) = V(E_{2r}(\Psi)) = V(\Pf_{2r}(\Psi)).
\end{equation}
In other words, for $n$ even, we have 
$\RR^1_{2d+1}(\Psi)=\RR^1_{2d}(\Psi)= \PP_{d}(\Psi)$, 
while for $n$ odd, we have
$\RR^1_{2d}(\Psi)=\RR^1_{2d-1}(\Psi)= \PP_{d-1}(\Psi)$. 

For $n=2k$, the Pfaffian of the generic skew-symmetric 
matrix $\Phi_n$ is given by 
\begin{equation}
\label{eq:pfaff}
\Pfaff(\Phi_n) = \sum_m \sigma(m) \omega(m),
\end{equation}
with the sum over all perfect matchings
$m=\{\{i_1,j_1\},\{i_2,j_2\},\dots,\{i_k,j_k\}\}$, 
(partitions of
$[2k]$ into blocks of size two with $i_p<j_p$), 
and where $\sigma(m)$ is the sign
of the corresponding permutation 
$\big(\begin{smallmatrix}1 & 2 & 3 & 4 &
\dots & 2k-1 & 2k\\ i_1 & j_1 & i_2 & j_2 & \dots & i_k &j_k 
\end{smallmatrix}\big)$, 
and $\omega(m) = y_{i_1,j_1} y_{i_2,j_2} \cdots
y_{i_k,j_k}$, see for instance \cite{BR}.  Note that 
$\Pfaff(\Phi_n)$ is a polynomial of degree $k=n/2$ in the variables $y_{i,j}$.  

For arbitrary $n$, it is known that the Pfaffian variety $\PP_d(\Phi_n)$ 
is irreducible, with singular locus $\PP_{d+1}(\Phi_n)$, see \cite{BE,KL}.  
These facts, together with Proposition \ref{prop:res bdry generic} 
and \eqref{eq:pfaff res}, yield the following.

\begin{corollary}
\label{cor:irr res}
Let $M_n$ be the boundary manifold of the 
general position arrangement $\A_n$.  For $n\ge 4$ 
and $\binom{n}{2} < d < \binom{n}{2}+n-2$, 
the resonance variety $\RR_d^1(M_n)$ is a singular, 
irreducible variety.
\end{corollary}

\begin{ack}
A significant portion of this work was done while the authors 
attended the Fall 2004 program ``Hyperplane Arrangements and 
Applications'' at the Mathematical Sciences Research Institute 
in Berkeley, California.  We thank MSRI for its support and 
hospitality during this stay.
\end{ack}

\end{document}